\documentclass[12pt]{article}

\usepackage{cite, epic,eepic,euscript, verbatim,amsmath,amssymb,amsthm,amsfonts,afterpage, float, bm,stmaryrd,paralist, bm,mathrsfs, color, multirow, cases}
\usepackage{booktabs,enumitem}
\usepackage{csvsimple}
\usepackage{colordvi}
\usepackage{graphicx}
\usepackage{subfigure}
\usepackage{caption}
\usepackage{wrapfig}
\usepackage{algpseudocode}
\usepackage{authblk}
\usepackage{lipsum}
\usepackage[titletoc]{appendix}
\usepackage{setspace}
\usepackage{lineno}
\usepackage{ulem}
\usepackage[size=tiny]{todonotes}
\usepackage[colorlinks=true,urlcolor=blue,allcolors=black]{hyperref}
\newtheorem{theorem}{Theorem}[section]
\newtheorem{lemma}{Lemma}[section]
\newtheorem{proposition}{Proposition}[section]
\newtheorem{corollary}{Corollary}[section]
\newtheorem{definition}{Definition}[section]
\newtheorem{remark}{Remark}[section]

\numberwithin{equation}{section}

\def\XXint#1#2#3{{\setbox0=\hbox{$#1{#2#3}{\int}$}
\vcenter{\hbox{$#2#3$}}\kern-.51\wd0}}

\newcommand{\R}{\mathbb{R}}

\newcommand{\ve}{\varepsilon}
\newcommand{\cA}{\mathcal{A}}
\newcommand{\cE}{\mathcal{E}}

\newcommand{\cL}{\mathcal{L}}
\newcommand{\cT}{\mathcal{T}}
\newcommand{\wc}{\rightharpoonup}

\newcommand{\1}{\mathbf 1}
\newcommand{\dist}{\operatorname{dist}}

\newcommand{\supp}{\operatorname{supp}}

\begin{document}

\title{Maximizing Reaction via Spatial Localization: 
The Reaction-Flux Landscape and Maximizing Sequences}

\author{Mohammad Yasir Feroz Khan}
\author{Bo Li}

\affil{Department of Mathematics\\
University of California San Diego\\
La Jolla, California, USA
\\
Email: m5khan@ucsd.edu and bli@math.ucsd.edu
}

\date{ }
\maketitle

\begin{abstract}
We continue our study of a PDE-constrained variational model for maximizing the total
reaction flux by optimizing the spatial distribution of catalytic materials, such as enzymes. The total amount of enzyme is fixed, and the substrate concentration is governed by a reaction–diffusion equation with prescribed boundary values. We consider both linear and nonlinear Michaelis–Menten reaction kinetics. In previous work, we constructed a maximizing sequence, 
computed the first and second variations of the reaction functional, and established the nonexistence of local maximizers. Here, we conduct a more comprehensive analysis. We prove that the functional is strictly concave for linear kinetics but need not be concave or convex for nonlinear kinetics. We also show that no critical points exist. We then characterize maximizing sequences and identify all their possible limiting measures. Our analysis reveals the mechanism of reaction maximization through spatial localization: maximizing sequences concentrate enzymes near the portion of the boundary where the substrate concentration is largest. Moreover, the landscape of the functional defining the total reaction differs between the linear and nonlinear cases, a point less emphasized in applications.
Consistent with existing biophysical studies, these findings motivate the development and analysis of regularized models.

\vspace{2 mm}

\vspace{2 mm}
\end{abstract}

\allowdisplaybreaks

\section{Introduction}
\label{s:Intro}

Chemical reactions, in which catalysts such as enzymes 
facilitate the conversion of substrates into products, are ubiquitous in chemistry and biology. In well-mixed systems, the efficiency of these processes is often measured by reaction rates 
\cite{Cornish-Bowden_2012,Lente_2015}.
In spatially heterogeneous systems such as biological cells, however, enzyme distributions can significantly 
influence reaction efficiency through the interplay of substrate diffusion, reaction, and loss. Experimental and theoretical biophysical studies have shown that spatial arrangements of enzymes in cells, 
including localization and co-localization, clustering and co-clustering, sequestration, and compartmentalization, 
can optimize the reaction flux under various biological constraints 
\cite{Gerland_CommunPhys2022,Gerland_JCP2013,Gerland_PRL2013,Gerland_NatPhys2022,Wingreen_NatBiotech2014,Hinzpeter_BiophysJ2017,tenWolde_PoLS2007,An_Biochem2017,Gerland_Interface2019}.
Optimizing enzyme distributions is particularly relevant to cell signaling pathways, which involve sequences of enzymatic reactions essential to biological functions
\cite{Lente_2015,Marcotte_Rev2012,Srere_EnzymeComplexes_Rev1987,Chenoweth_JACS2018}.

We study the variational problem of maximizing the reaction functional
of the total reaction flux
\begin{equation}
\label{Ju}
J[u]=\int_\Omega u f(\rho_u)\,dx
\end{equation}
over the class of admissible enzyme concentrations
\begin{equation}
\label{cE}
\cE=\left\{u\in L^2(\Omega):u\ge0\ \mbox{a.e.\ in }\Omega
\ \mbox{and}\ \int_\Omega u\,dx=M_{\rm e}\right\},
\end{equation}
where $\Omega\subset\R^3$ is a bounded domain 
and $M_{\rm e}>0$ is a constant representing the prescribed total amount of enzyme. Under the assumptions stated in section~\ref{s:Model}, for each $u\in\cE$, the substrate concentration $\rho_u$
is the unique nonnegative weak solution to the boundary-value problem 
\begin{equation}
\label{BVP4rho}
\left\{
\begin{aligned}
-\Delta\rho_u+\sigma\rho_u+u f(\rho_u)&=0 &&\quad \mbox{in }\Omega,\\
\rho_u&=g&&\quad \mbox{on }\partial\Omega.
\end{aligned}
\right.
\end{equation}
Here $\sigma>0$ is a constant and $g \ge 0$ is a given continuous function on the boundary $\partial \Omega$.
The diffusion and reaction coefficients have been rescaled to one.
The constant $\sigma$ describes substrate loss other than the
reaction modeled by $uf(\rho_u)$.
We consider two commonly used reaction kinetics given in the dimensionless form
\begin{equation}
\label{f}
f(s)=\left\{
\begin{aligned}
&s&&\quad \mbox{for the linear case},\\
&\frac{s}{1+s}&&\quad \mbox{for the nonlinear case}.
\end{aligned}\right.
\end{equation}
The nonlinear case is the 
Michaelis--Menten reaction kinetics.

This PDE-constrained variational model was constructed in our previous work \cite{KhanLiFirst}, following
the biophysical studies
\cite{Gerland_PRL2013,Gerland_JCP2013} (cf.\ also
\cite{Gerland_CommunPhys2022,Gerland_NatPhys2022,Rios_JCP2025,Plasson_PRE2025}). 
As in \cite{Gerland_PRL2013,Gerland_JCP2013}, we consider the steady-state substrate
concentration $\rho_u$; cf.\ \eqref{BVP4rho}. 
Focusing on the localization, we use the Dirichlet instead of Neumann boundary conditions. 
While other effects, such as transport via fluid 
flow, electrostatic interaction, 
and a multi-stage reaction process, are excluded, our analysis will show that 
the PDE-constrained variational model has the essential components
of reaction to facilitate the catalyst spatial localization that can be experimentally implemented
\cite{Chenoweth_JACS2018}.

The main results we obtained in our previous work \cite{KhanLiFirst} include the following: 
\begin{compactenum}
    \item[(1)]
    The well-posedness of the boundary-value problem \eqref{BVP4rho} and the continuous dependence of
    the reaction functional $J[u]$ on the enzyme concentration $u$; 
    \item[(2)] 
    The reaction functional $J: \cE \to \R$ has the supremum $M_{\rm e} f(g_{\rm max})$, which is proved by
    the construction of a maximizing sequence, but it is not attainable
   \begin{equation}
\label{supJ}
M_{\rm e}f(g_{\rm max}) = \sup_{v \in \cE} J[v] > J[u] \qquad\forall u \in \cE, 
\end{equation}
where $ g_{\rm max}=\max_{\partial\Omega}g$, assumed to be positive. 
\item[(3)]
Explicit formulas of the first and second variations of the reaction functional $J: \cE\to \R;$
\item[(4)]
The reaction functional $J: \cE \to \R$ permits no local maximizers with respect to the $L^2(\Omega)$-norm. 
\end{compactenum}
Some of these results will be recalled in section~\ref{s:Model}.

Here in this work, we present a more comprehensive and relatively complete 
analysis of the PDE-constrained variational problem and interpret our analysis in terms of 
the mechanisms underlying the reaction maximization. 
Our main results and the key ideas of their proofs are outlined below.

\uline{(1)  Concavity 
of 
the reaction functional and the difference between the linear and nonlinear kinetics.}
We prove that the reaction functional $J: \cE \to \R$ is strictly concave for the linear kinetics; cf.~Theorem~\ref{t:JConcavity}. We also construct an example 
of a system with special $\Omega$, $M_{\rm e}$, and $g$ 
to show that
the functional $J: \cE \to \R$ is neither concave nor convex for the nonlinear kinetics; cf.\ Theorem~\ref{t:MM}. 
These results are obtained using the explicit formulas of the first and second variations of the reaction 
functional. The linearity simplifies the formula while the nonlinearity requires 
more delicate and tedious calculations. 

The resulting structural difference between the linear and nonlinear kinetics, often neglected in application,
may have crucial consequences in modeling a reaction system. For linear kinetics, adding the gradient term
$-(\ve/2)|\nabla u|^2$ to the integrand defining $J[u]$,
with $\ve>0$, gives a strictly concave functional on
$H^1(\Omega)\cap\cE$. Existence of a global maximizer
follows by the direct method, and strict concavity
gives uniqueness.
For the nonlinear
kinetics, the example shows that the regularized functional with a small $\ve  > 0$ may still be 
neither concave nor convex. If so, such a nonlinear model will be more useful in application, as it can allow
the transition between different stable states once fluctuations are included further in the model.

\begin{wrapfigure}[11]{R}{46 mm}
\vspace{-12 mm}

\includegraphics[width=45mm,height=45mm]{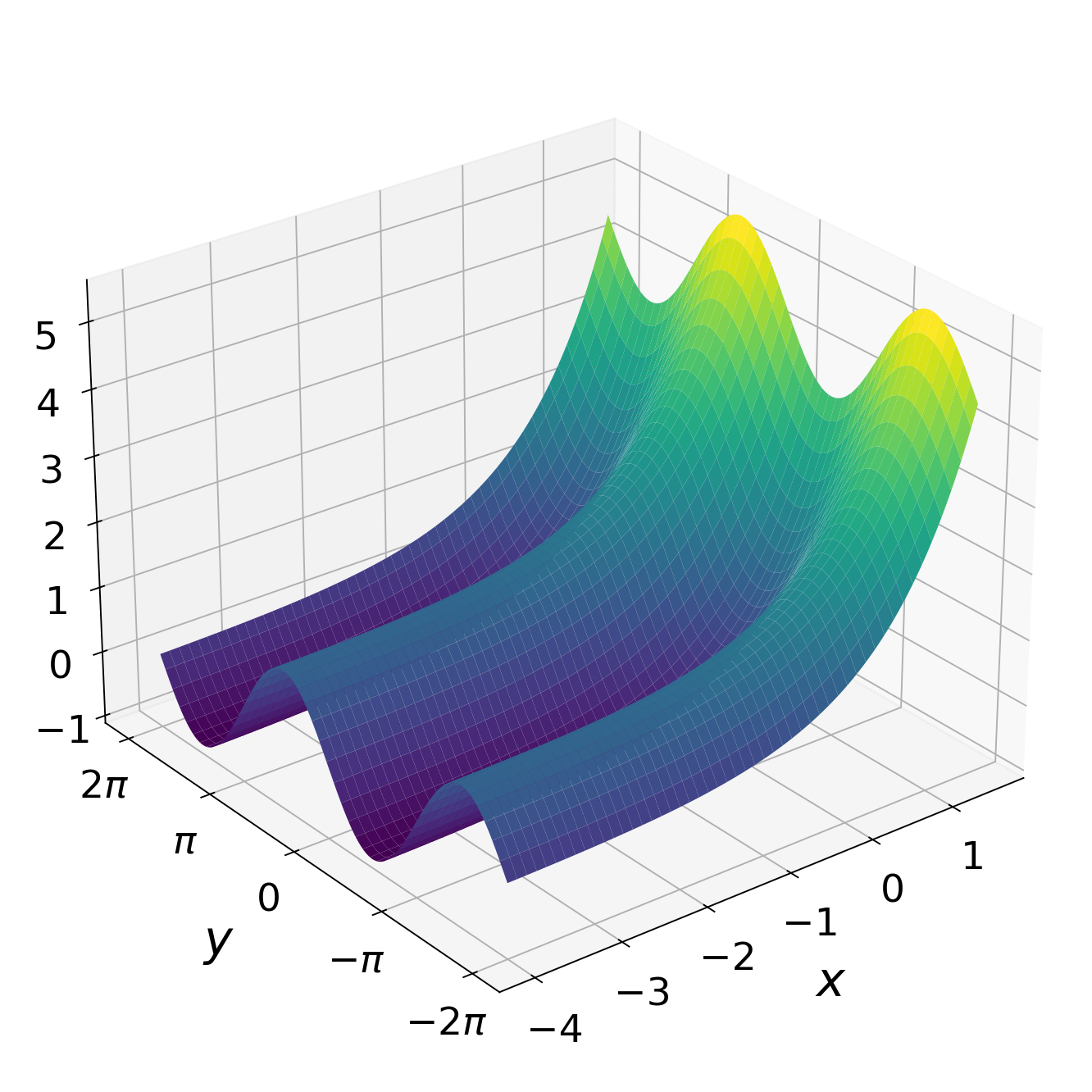}
\vspace{-4 mm}
\caption{\footnotesize{The function $h(x,y)= e^x+\sin y$ is neither
concave nor convex and has no critical points.}}
\label{f:landscape}
\end{wrapfigure}
$\, \! \! \!$
\uline{(2) Non-existence of critical points.}
We prove that the reaction functional $J:\cE\to\R$
admits no critical points for both the linear and nonlinear kinetics; cf.\ Theorem~\ref{t:CriticalPoint}.
More precisely, we prove that for any $u \in \cE$, there exists an admissible perturbation $v_u \in L^2(\Omega)$
such that $J'[u][v_u] > 0.$ The proof uses the explicit formula of the first variation $J'[u]$ 
and the maximum principle applied to an auxiliary function.

The non-existence of a critical point together with the non-concavity and non-convexity
of the reaction functional $J: \cE \to \R$ in the nonlinear example of Theorem~\ref{t:MM} imply that
the energy landscape has some special features. 
For instance, it may be oscillatory in some directions
while monotonic in other directions.  Figure~\ref{f:landscape} shows the landscape of
a finite-dimensional "energy functional" which has such features. For infinite-dimensional 
functionals, the landscape can be possibly more complicated.

\uline{(3) Characterization of maximizing sequences.}
We prove that, if $\{ u_k \} \subset \cE $  is a maximizing sequence of the functional $J: \cE \to \R$, 
then $\| u_k \|_{L^2(\Omega)} \to \infty$ and $\rho_{u_k} \to \rho_0$ in $H^1(\Omega)$, 
where $\rho_{u_k}$ and $\rho_0$ are the substrate concentrations defined by \eqref{BVP4rho} with $u=u_k$ and $u = 0,$ respectively, and
\begin{equation}
\label{Kuk}
\int_{K\cap \Omega} u_k \, dx \to 0 
\end{equation}
for any compact subset $K \subset \overline{\Omega} \setminus \Gamma_{\rm max}$, where
$
\Gamma_{\rm max}=\{y\in\partial\Omega:g(y)=g_{\rm max}\}. 
$
Thus, to maximize the reaction flux, it is necessary to distribute enzymes near $\Gamma_{\rm \max}$ where substrate supply is the largest.

Furthermore,  we 
prove that $\{u_k\}\subset\cE$ is a maximizing sequence 
if and only if 
\[
\int_\Omega\bigl[ f(g_{\rm max})-f(\rho_0)\bigr] u_k\,dx\to0
\]
and there exist $0\le v_k\le u_k$ $(k \ge 1)$ with
$\|u_k-v_k\|_{L^1(\Omega)}\to0$ and
$\|v_k\|_{H^{-1}(\Omega)}\to0$;
cf.\ Theorem~\ref{t:ukJ}.  By choosing $v_k = u_k$ $(k \ge 1)$, we obtain sufficient conditions.
The canonical maximizing sequence $\{ u_k \} \subset \cE$
we constructed in our previous work \cite{KhanLiFirst} satisfies such sufficient conditions. However, 
this is not true in general; cf.\ Theorem~\ref{t:sufficient}.

\begin{figure}[!htbp]
\centering
\includegraphics[width=0.8\textwidth]{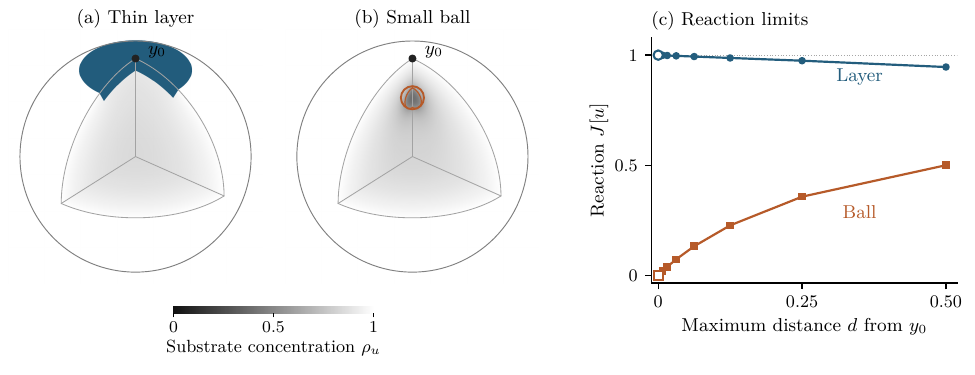}
\vspace{-4 mm}

\caption{\small{Different limits of total reaction flux for the same 
limiting enzyme concentration $M_{\rm e} \delta_{y_0}$ with $y_0 \in \Gamma_{\rm max}.$ 
The colored regions show enzyme concentrated in
(a) a thin region near the boundary and
(b) a small interior ball. Both sequences have the measure limit
$M_{\rm e}\delta_{y_0}$.
(c) Numerical computations using the unit ball $\Omega \subset \R^3$, 
$f(s) = s$, and $M_{\rm e}=\sigma=g=1$ show that, as 
the maximum distance $d$ of the enzyme support
from $y_0$ tends to zero, the total reaction flux approaches
$M_{\rm e}f(g_{\rm max}) = \sup_{\cE} J$ for the thin regions and zero
for the small balls.}}
\label{f:concentration}
\end{figure}

\uline{(4) The limiting Radon measures.} By \eqref{Kuk}, we expect every weak-$\ast$ subsequential limit of the measures associated with a maximizing sequence $\{u_k\}$ to be a nonnegative Radon measure of mass $M_{\rm e}$ supported on $\Gamma_{\rm \max}$. This is indeed true. Conversely, any such measure can be shown to be the weak-$\ast$ limit of a maximizing sequence.
More generally, for any such measure $\nu$ and any $\ell$ such that 
$0 \le \ell \le M_{\rm e}f(g_{\rm max}) = \sup_{\cE} J$, there exists a sequence $\{ v_k \} \subset \cE$ such
that $v_k\,dx\stackrel{*}{\rightharpoonup}d\nu$ on
$\overline\Omega$ 
with $J[v_k] \to \ell;$ cf.\ Theorem~\ref{t:RadonMeasures}. So, 
the geometry of enzyme concentration affects substrate depletion
even when the limiting enzyme measure is fixed.
The proof is done by first constructing two sequences in $\cE$ that both weakly-$\ast$ converge to the same
measure but with the functional sequences converging to $0$ and $M_{\rm e}f(g_{\rm max})$, 
using small balls inside $\Omega$ and thin layers near $\partial \Omega$, respectively, 
cf.\ Figure~\ref{f:concentration}. The general case $\ell$ is proved by interpolation.

\medskip

Our analysis suggests constructing a regularized functional that permits a global maximizer. In our previous work, 
we proposed such a regularized functional
\[
J_{\ve,\beta}[u] = J[u] -\frac{\ve}{2 } \int_\Omega |\nabla u|^2 \, dx -\frac{1}{\beta}\int_\Omega u \log u \, dx,
\]
where both $\ve > 0$ and $\beta > 0$ are constants, and $u \in \cE \cap H^1(\Omega).$ 
It is shown using the direct method in the calculus of variations that the regularized functional has a global maximizer \cite{KhanLiFirst}.
We will present more detailed studies in our upcoming work.

The rest of the paper is organized as follows. In
section~\ref{s:Model}, we recall the results of
\cite{KhanLiFirst} that will be used and prove additional results. 
In section~\ref{s:Concavity}, we prove the strict concavity for the linear
kinetics and construct an example to show the non-concavity and non-convexity for the nonlinear kinetics. 
In section~\ref{s:CriticalPoints},
we prove the non-existence of critical points.
In section~\ref{s:Characterization}, we characterize
maximizing sequences of the functional $J$. Finally, in section~\ref{s:LimitMeasures}, we study the limiting
Radon measures.

\section{Some Preliminary Results}
\label{s:Model}

In this section, we first make assumptions and set notations, then review the main results 
of our previous work \cite{KhanLiFirst} that will be used here in our continued analysis, and finally
add some new and basic results on the reaction functional.

We assume the following throughout:
\begin{compactenum}
\item[(A1)]
The set $\Omega\subset\R^3$ is a bounded domain (i.e., an open and connected set) with a
$C^2$-boundary $\partial \Omega.$

\item[(A2)]
Both $M_{\rm e}>0$ and $\sigma>0$ are given constants.

\item[(A3)]
The function $g\in H^2(\Omega)\hookrightarrow C(\overline\Omega)$
is nonnegative on $\overline\Omega$, and
\begin{equation}
\label{gmax}
g_{\rm max}:=\max_{\partial\Omega}g
=\|g\|_{L^\infty(\partial\Omega)}>0.
\end{equation}
\end{compactenum}

For any $p \in [1, \infty]$ and any $h\in H^1(\Omega)$, we denote
\begin{align*}
&L^p_+ (\Omega) = \{ u \in L^p(\Omega): u \ge 0 \mbox{ a.e.\ in } \Omega \}, 
\\
&
H^1_h(\Omega) = \{ v  \in H^1(\Omega): v= h \mbox{ on } \partial \Omega \}, 
\end{align*}
where the boundary value is understood as trace. For any
Lebesgue-measurable set $S\subseteq\R^3$, we denote by $|S|$ its
Lebesgue measure. 
We shall use \eqref{Ju} to define $J[u]$ for any
$u\in L^2_+(\Omega)$, 
denote by $\rightharpoonup$ the weak convergence with respect to an underlying topology, and 
denote by $\log $  the natural logarithm with $\log s = -\infty$ if $s \le 0$ and
$s \log s = 0$ if $s = 0.$

The following proposition summarizes some basic properties of the 
solution to the boundary-value problem \eqref{BVP4rho} and 
the reaction functional $J: \cE\to \R:$

\begin{proposition}
\label{p:FirstPaper}
The following statements hold:
\begin{compactenum}
\item[{\rm (1)}]
Let $u\in L^2_+(\Omega)$. The boundary-value problem
\eqref{BVP4rho} has a unique nonnegative weak solution, denoted
$
\rho_u=\cT_g(u)\in H_g^1(\Omega)\cap H^2(\Omega).
$
Moreover, $0<\rho_u<g_{\rm max}$ in $\Omega$.

\item[{\rm (2)}]
If $u_k,u\in L^2_+(\Omega)$ and $u_k\rightharpoonup u$ in
$L^2(\Omega)$, then
$\rho_{u_k}\to\rho_u$ in $H^1(\Omega)$ and $J[u_k]\to J[u].$

\item[{\rm (3)}]
Let $u,v\in L^2_+(\Omega)$ be such that $u\le v$ a.e.\ in
$\Omega$. Then
$\rho_u\ge\rho_v$ in $\Omega$, $\partial_\nu \rho_u \le \partial_\nu\rho_v$ a.e.\
on $\partial \Omega$,
where $\nu$ is the unit normal at $\partial \Omega$ exterior to 
$\Omega$,  and $J[u] \le J[v].$ 
\end{compactenum}
\end{proposition}

\begin{proof}
(1) The existence, uniqueness, and regularity of the solution
$\rho_u=\cT_g(u)$ are proved in \cite{KhanLiFirst}
using the direct method in the calculus of variations and
elliptic regularity \cite{GilbargTrudinger98}.

We can rewrite the equation as
\[
-\Delta\rho_u+ \left[ \sigma+a(u) \right] \rho_u=0 \qquad \mbox{in } \Omega, 
\]
where $a(u) = u$ for the linear case and $a(u) = u / ( 1 + \rho_u)$ for the nonlinear case, 
respectively. Note that $d := \sigma + a(u) \ge 0$ a.e.\ in $\Omega$ and $ d \in L^2(\Omega).$ Since $\rho_u\ge0$, the equation gives
$\Delta\rho_u=d\rho_u\ge0$.
The maximum principle \cite{GilbargTrudinger98} immediately implies that $0 \le \rho_u \le g_{\rm max}$ on $\Omega.$

The strict bound $\rho_u < g_{\rm max}$ 
follows from the maximum principle
\cite{Trudinger1973} (Corollary~5.1), and the fact that $\rho_u $ is not a constant, for otherwise
it would be $0$ in $\overline{\Omega}$ by \eqref{BVP4rho}, contradicting 
$\rho_u = g$ on $\partial \Omega$ with $g_{\rm max} > 0.$
We claim that $\rho_u>0$ in $\Omega$. Suppose there exists $x_0 \in \Omega$ such that $\rho_u(x_0)=0.$ 
Then, the set 
$Z: =\{x\in\Omega:\rho_u(x)=0\}$ is not empty. 
Since $\rho_u$ is continuous, $Z$ is closed in $\Omega$.
For any $x\in Z$, choose $r>0$ such that
$B_{5r}(x)\subset \subset \Omega$. Since
$\sigma+a(u)\in L^2(\Omega)$, the local Harnack inequality
(cf.\ Corollary~5.2 in \cite{Trudinger1973}) gives
\[
0\le\sup_{B_r(x)}\rho_u
\le C\inf_{B_r(x)}\rho_u=0.
\]
Hence, $\rho_u \equiv0$ on $B_r(x)$. Thus, $B_r(x)\subset Z$. Since $x\in Z$ is arbitrary,
$Z$ is also open in $\Omega$. By the connectedness of
$\Omega$, $Z=\Omega$. Thus, $\rho_u\equiv0$ in $\Omega$, and on $\overline\Omega$ by continuity,
contradicting $\rho_u=g$ on $\partial\Omega$ and
$g_{\rm max}>0$. Hence $\rho_u>0$ in $\Omega$.

(2) The weak convergence $u_k\rightharpoonup u$ in $L^2(\Omega)$ implies
$\sup_k\|u_k\|_{L^2(\Omega)}<\infty$. Note that $ f(\rho_{u_k})\le f(g_{\rm max})$
in $\Omega$.
By the global estimate for the Laplacian with the Dirichlet boundary condition
(cf.\ Theorem~8.12 in \cite{GilbargTrudinger98}), we have for all $k \ge 1$ that
\begin{align*}
\|\rho_{u_k}\|_{H^2(\Omega)}
&\le C\left(
\|\rho_{u_k}\|_{L^2(\Omega)}
+\|\sigma\rho_{u_k}+u_kf(\rho_{u_k})\|_{L^2(\Omega)} + \|g\|_{H^2(\Omega)}
\right)
\le C,
\end{align*}
for some generic constant $C > 0.$
By the compact embeddings of $H^2(\Omega)$ into
$H^1(\Omega)$ and $C(\overline\Omega)$, there exist
a subsequence, not relabeled, and $\rho\in H^2(\Omega)$
such that
$
\rho_{u_k}\rightharpoonup\rho $ in $ H^2(\Omega)$ and 
$\rho_{u_k}\to\rho $ in $H^1(\Omega)$ and $C(\overline\Omega).$
In particular, $0\le\rho\le g_{\rm max}$ on $\Omega$
and $\rho=g$ on $\partial\Omega$, so $\rho\in H_g^1(\Omega)$.
By the definition of $\rho_{u_k} = \cT_g(u_k)$, we have for 
 any $k \ge 1$ that 
\begin{equation}
\label{rhouk}
\int_\Omega
\left[
\nabla\rho_{u_k}\cdot\nabla\varphi
+\sigma\rho_{u_k}\varphi
+u_k f(\rho_{u_k})\varphi
\right]dx=0 \qquad \forall \varphi\in H_0^1(\Omega).
\end{equation}
Sending $k\to \infty$ and using $f(\rho_{u_k})\to f(\rho)$ in $L^\infty(\Omega)$, we obtain
\[
\int_\Omega
\left[
\nabla\rho\cdot\nabla\varphi
+\sigma\rho\varphi
+u f(\rho)\varphi
\right]dx=0
\qquad\forall\varphi\in H_0^1(\Omega).
\]
By uniqueness in part~(1), $\rho=\rho_u$.
It now follows from $ \sup_{k\ge 1} \| u_k \|_{L^2(\Omega)} < \infty$, 
$\rho_{u_k} \to \rho_u$ in $C(\overline{\Omega})$, and 
$u_k \rightharpoonup u$ in $L^2(\Omega)$ that
\[
J[u_k]-J[u] = \int_\Omega u_k \bigl[f(\rho_{u_k})-f(\rho_u)\bigr]\,dx
+\int_\Omega(u_k-u)f(\rho_u)\,dx \to 0 \qquad \mbox{as } k \to \infty.
\]
Every subsequence admits a further subsequence with the desired convergence. 
Hence, the same convergence holds for the full sequence. 

(3) This is Theorem~2.2  in \cite{KhanLiFirst}. 
\end{proof}

Let $u\in L^2_+(\Omega)$. We introduce the linear and self-adjoint operator \cite{KhanLiFirst}
\begin{equation}  
\label{cLu}
\cL_u=(-\Delta+\sigma+u f'(\rho_u))^{-1}:
L^2(\Omega)\longrightarrow H_0^1(\Omega)\cap H^2(\Omega).
\end{equation}
Note that the $H^2(\Omega)$-regularity follows from a bootstrapping argument \cite{KhanLiFirst}. 
The operator $\cL_u$ is bijective. 
The function $\psi_u := \cL_u \bigl(u f'(\rho_u)\bigr)\in H_0^1(\Omega)\cap H^2(\Omega) $
is the unique solution to the boundary-value problem
\begin{equation}
\label{BVPpsi}
\left\{
\begin{aligned}
-\Delta\psi_u+(\sigma+u f'(\rho_u))\psi_u
&=u f'(\rho_u)&&\mbox{in }\Omega,\\
\psi_u&=0&&\mbox{on }\partial\Omega.
\end{aligned}
\right.
\end{equation}

\begin{proposition} 
\label{p:psiu}
Let $u \in L_+^2(\Omega)$, $\rho_u = \cT_g(u) \in H_g^1(\Omega) \cap H^2(\Omega)$, 
and $\psi_u: = \cL_u \bigl( u f'(\rho_u) \bigr) \in H_0^1(\Omega) \cap H^2(\Omega).$ Then, we have 
$0 \le \psi_u < 1$ on  $\overline{\Omega}.$
\end{proposition}

\begin{proof}
By the direct method in the calculus of variations \cite{EvansBook2010,GilbargTrudinger98}, we verify that 
$\psi_u $ is the unique minimizer of $K_u: H_0^1(\Omega) \to \R\cup \{ \infty \}$, defined by
\[
K_u[\eta] = \int_\Omega \left[\frac12 |\nabla \eta|^2 
+ \frac{\sigma}{2} \eta^2 + \frac12 u f'(\rho_u) \eta^2 - u f'(\rho_u) \eta \right]dx 
\qquad \forall \eta \in H_0^1(\Omega); 
\]
see \cite{KhanLiFirst} for more details. 
Setting $\psi = \max \{ 0, \min \{ \psi_u, 1 \} \} \in H_0^1(\Omega)$, we obtain $ 0 \le \psi \le 1$ 
on $\overline{\Omega}$, and
\begin{align*}
0 & \ge K_u[\psi_u] - K_u [\psi] \nonumber \\
&   \ge \int_{ \{ \psi_u <  0 \}} 
\left[ \frac{\sigma}{2} \psi_u^2  + \frac{1}{2} u f'(\rho_u) \psi_u (\psi_u - 2) \right] dx
\\
&\qquad + \int_{ \{ \psi_u > 1\}} 
\left[ \frac{\sigma}{2} \left( \psi_u^2 - 1\right) 
+ \frac{1}{2} u f'(\rho_u) \left( \psi_u -1 \right)^2  \right]dx 
\nonumber \\
&\ge 0. 
\end{align*}
Thus, $K_u[\psi_u] = K_u [\psi]$ and hence  $\psi$ is also a minimizer of $K_u$ over $H_0^1(\Omega)$. The
uniqueness of such a minimizer implies that $\psi_u = \psi \in [0, 1]$ on $\overline{\Omega}$.

To prove the strict upper bound, we set $z=1-\psi_u-\rho_u/g_{\rm max} \in H^2(\Omega)$. By
\eqref{BVP4rho} and \eqref{BVPpsi}, we have
\[
\left\{
\begin{aligned}
-\Delta z+(\sigma+uf'(\rho_u))z
&=\sigma+\frac{u}{g_{\rm max}}
\bigl[f(\rho_u)-\rho_u f'(\rho_u)\bigr]
&&\mbox{ in }\Omega,\\
z&=1-\frac{g}{g_{\rm max}}\ge0
&&\mbox{ on }\partial\Omega.
\end{aligned}
\right.
\]
For both reaction laws in \eqref{f},
$f(s)-sf'(s)\ge0$ for $s\ge0$.
Thus, the right-hand side of the above equation in $\Omega$ is nonnegative a.e.\ in $\Omega.$
It is also in $L^2(\Omega).$
The coefficient of $z$ in the left-hand side of the equation is strictly positive and
is in $L^2(\Omega)$. Therefore,
by the weak minimum principle
(cf.\ Corollary~3.1 in \cite{Trudinger1973}), we have $z\ge0$ in $\Omega$. Since $\rho_u>0$ in $\Omega$ by
Proposition~\ref{p:FirstPaper}, we obtain
$\psi_u \le 1 - \rho_u / g_{\rm max} < 1$ in $\Omega.$
Together with  $\psi_u=0$ on
$\partial\Omega$, this proves that $\psi_u < 1 $ on $\overline{\Omega}.$
\end{proof}

Let  $u\in\cE$ (cf.\ \eqref{cE} for the definition of $\cE$). We define the set $\cA_u$
of admissible perturbations for $u$ by
\begin{equation}
\label{cAu}
\cA_u=\left\{v\in L^2(\Omega):
u+tv\in\cE\ \mbox{if }|t|\le t_v
\ \mbox{for some }t_v>0\right\}.
\end{equation}
The following proposition is proved in Theorems~4.1--4.3 of \cite{KhanLiFirst}.

\begin{proposition}
\label{p:variations} Let $u \in \cE$, 
$\rho_u = \cT_g(u) \in H_g^1(\Omega) \cap H^2(\Omega)$, and
 $\psi_u: = \cL_u \bigl( u f'(\rho_u) \bigr) \in H_0^1(\Omega) \cap H^2(\Omega).$
For any $v \in \cA_u$, let $\hat{\psi}_{u,v}=\cL_u\bigl(-v f(\rho_u)\bigr).$ 
Then, we have the first and second variations
\begin{align}
\label{J'uv}
&J'[u][v]=\int_\Omega f(\rho_u)(1-\psi_u)v\,dx
\qquad\forall v\in\cA_u, 
\\
\label{J''uv}
&J''[u][v]
=\int_\Omega\left[
2v f'(\rho_u)(1-\psi_u)\hat{\psi}_{u,v}
+u f''(\rho_u)(1-\psi_u)\hat{\psi}_{u,v}^2
\right]dx
\qquad\forall v\in\cA_u.
\end{align}
Let $\alpha_v(t) = J[u+tv]$ $(|t| < t_v).$ Then, $\alpha_v \in 
C^2((-t_v, t_v))$. Moreover, 
$\alpha_v'(t) = J'[u+tv][v] $ and $\alpha_v''(t) = J''[u+tv][v]$  for any $t \in (-t_v, t_v)$, and 
\begin{equation*}
J[u+tv] = J[u] + J'[u][v]t + \frac12 J''[u][v] t^2 + o(t^2) \qquad \mbox{as } t \to 0.
\tag*{\qed}
\end{equation*}
\end{proposition}

We finally prove the Lipschitz-continuity of the functional $J: \cE \to \R$ with respect to 
the $L^1(\Omega)$-norm. 

\begin{proposition}
\label{p:Juv}
We have
\begin{equation*}
|J[v]-J[u]|
\le f(g_{\rm max}) \|v-u\|_{L^1(\Omega)}
\qquad\forall u,v\in\cE.
\end{equation*}
\end{proposition}

\begin{proof} 
Let $v, u \in \mathcal E.$
Set $u_t=(1-t)u+tv\in\cE$ for $0\le t\le1$.
Fix $0<t<1$. Then, 
\[
u_t+s(v-u)=(1-t-s)u+(t+s)v\in\cE \qquad \mbox{if } |s|<\min\{t,1-t\}. 
\]
Thus, $v-u \in \cA_{u_t}$. 
 By Proposition~\ref{p:FirstPaper}, Proposition~\ref{p:psiu}, and the monotonicity of $f$, 
 \[
0\le f(\rho_{u_t})(1-\psi_{u_t})
\le f(g_{\rm max})\qquad \mbox{in } \Omega. 
 \]
Consequently, 
\begin{align}
\left| \frac{d}{dt}J[u_t] \right|
& = \left| J'[u_t][v-u] \right| \nonumber \\
&= \left| \int_\Omega f(\rho_{u_t})(1-\psi_{u_t})(v-u)\,dx \right| \nonumber \\
& \le f(g_{\rm max}) \| u - v \|_{L^1(\Omega)}
\qquad \forall t \in (0, 1). 
\label{JJ'ut}
\end{align}
By Proposition~\ref{p:FirstPaper} and
\eqref{JJ'ut}, $t\mapsto J[u_t]$
is continuous on $[0,1]$ and differentiable on $(0,1)$.
The mean value theorem and \eqref{JJ'ut}
therefore imply the desired inequality. 
\end{proof}

\section{Concavity}
\label{s:Concavity}

In this section, we prove that the reaction functional $J: \cE \to \R $ is 
strictly concave for the linear case, i.e., $f(s) = s $ $(s \in \R)$; cf.\ Theorem~\ref{t:JConcavity}. We also
construct $\Omega$, $\sigma$, and $g$ of the model and show that the corresponding reaction functional 
$J: \cE \to \R$ 
is neither concave nor convex for the nonlinear case, i.e., $f(s) = s/(s+1)$ $(s > -1)$; cf.\ Theorem~\ref{t:MM}.

\begin{theorem}[Strict concavity of the reaction functional with the linear kinetics]
\label{t:JConcavity}
For $f(s)=s$ $(s \in \R)$, we have
\begin{equation}
\label{J''linear}
J''[u][v]<0\qquad\forall u\in\cE\quad
\forall v\in\cA_u\setminus\{0\}.
\end{equation}
Consequently, $J$ is strictly concave on $\cE$.
\end{theorem}

\begin{proof}
Let $u\in\cE$ and $v\in\cA_u\setminus\{0\}$. 
Let $\rho_u = \cT_g(u) \in H_g^1(\Omega)\cap H^2(\Omega)$ and 
$\psi_u = \cL_u(uf'(\rho_u)) = \cL_u(u) \in H_0^1(\Omega) \cap H^2(\Omega)$ be given 
by \eqref{BVP4rho} and \eqref{BVPpsi}, respectively. 
Denote for simplicity
\[
\rho=\rho_u, \qquad q=1-\psi_u,\qquad
\xi=\cL_u(vq),\qquad \theta=\frac\xi q, \qquad \hat{\psi}_{u,v}=\cL_u\bigl(-v f(\rho_u)\bigr), 
\]
where $\cL_u$ is defined in \eqref{cLu}. 
Here $\theta\in H_0^1(\Omega)\cap H^2(\Omega)$, since $q \in H^2(\Omega)$ and 
$ \min_{\overline{\Omega}} q >0$ by Proposition~\ref{p:psiu}. 
Since $f'=1$ and $f''=0$, it follows from the second variation formula in \eqref{J''uv} 
and the self-adjointness of $\cL_u$ that 
\begin{align}
\label{J''xi}
J''[u][v]
&=2\int_\Omega vq\hat{\psi}_{u,v}\,dx
=-2\int_\Omega vq\cL_u(v\rho)\,dx
\nonumber \\
&
=-2\int_\Omega v\rho\cL_u(vq)\,dx
=-2\int_\Omega v\rho\xi\,dx.
\end{align}

By \eqref{BVP4rho}, \eqref{BVPpsi}, and \eqref{cLu},
the equations satisfied by $\rho$, $q$, and $\xi$ are
\[
\begin{aligned}
-\Delta\rho+(\sigma+u)\rho&=0,\\
-\Delta q+(\sigma+u)q&=\sigma,\\
-\Delta\xi+(\sigma+u)\xi&=vq.
\end{aligned}
\]
Substituting $\xi=q\theta$ into the third equation and using the equation for $q$,
we obtain
\begin{align}
\label{BVPtheta}
vq
&=-\Delta(q\theta)+(\sigma+u)q\theta
\nonumber \\
&=-q\Delta\theta-2\nabla q\cdot\nabla\theta
-\theta\Delta q+(\sigma+u)q\theta
\nonumber \\
&=-q\Delta\theta-2\nabla q\cdot\nabla\theta
+\theta\bigl[-\Delta q+(\sigma+u)q\bigr]
\nonumber \\
&=\sigma\theta-q\Delta\theta
-2\nabla q\cdot\nabla\theta.
\end{align}
The equations for $\rho$ and $q$ also imply
\begin{align}
\label{divqrho}
\nabla \cdot (q\nabla\rho-\rho\nabla q)
&=q\Delta\rho-\rho\Delta q
\nonumber \\
&=q(\sigma+u)\rho
-\rho\bigl[(\sigma+u)q-\sigma\bigr]
\nonumber \\
&=\sigma\rho.
\end{align}
Multiplying \eqref{BVPtheta} by $ \rho\theta$ and integrating
by parts, we obtain by \eqref{divqrho} that
\[
\begin{aligned}
\int_\Omega v\rho\xi\,dx
&=\sigma\int_\Omega\rho\theta^2\,dx
-\int_\Omega\rho q\theta\Delta\theta\,dx
-2\int_\Omega\rho\theta\nabla q\cdot\nabla\theta\,dx\\
&=\sigma\int_\Omega\rho\theta^2\,dx
+\int_\Omega\rho q|\nabla\theta|^2\,dx
+\frac12 \int_\Omega \nabla \left( \theta^2 \right) \cdot (q\nabla\rho-\rho\nabla q) dx\\
&=\sigma\int_\Omega\rho\theta^2\,dx
+\int_\Omega\rho q|\nabla\theta|^2\,dx
-\frac12\int_\Omega
\theta^2 \nabla \cdot (q\nabla\rho-\rho\nabla q)\,dx\\
&=\int_\Omega\rho q|\nabla\theta|^2\,dx
+\frac\sigma2\int_\Omega\rho\theta^2\,dx. \qquad \mbox{[by \eqref{divqrho}]}
\end{aligned}
\]
Combining this with \eqref{J''xi}, we get
\begin{equation}
\label{J''theta}
J''[u][v]
=-2\int_\Omega\rho q|\nabla\theta|^2\,dx
-\sigma\int_\Omega\rho\theta^2\,dx\le0.
\end{equation}
If $J''[u][v]=0$, then $\rho,q>0$ in $\Omega$ and $\sigma>0$ would imply
$\theta=0$, and hence by \eqref{BVPtheta}
$vq=0$ a.e.\ in $\Omega$. 
Since $q>0$ in $\Omega$, this implies $v = 0$ a.e.\ in $\Omega$, contradicting the assumption that
$v\ne0$ in $L^2(\Omega)$. Thus \eqref{J''linear} holds.

Let $u_0,u_1\in\cE$, $u_0\ne u_1$, and set
\[
u_t=(1-t)u_0+tu_1,\qquad v=u_1-u_0,
\qquad \alpha(t)=J[u_t]\qquad(0\le t\le1).
\]
We have $v\in\cA_{u_t}$ if $0<t<1$.
By Proposition~\ref{p:FirstPaper}, Proposition~\ref{p:variations}, and Proposition~\ref{p:Juv},
$\alpha\in C([0,1])\cap C^2((0,1))$ and 
$\alpha''(t)=J''[u_t][v]<0 $ for all $t \in (0, 1).$
Therefore, $\alpha$ is strictly concave on $[0, 1]$, i.e., 
$\alpha(t) > (1-t) \alpha(0) + t \alpha(1) $ for all $t \in (0, 1)$. Consequently, 
\[
J[(1-t) u_0 + tu_1] >  (1-t) J[u_0] + t J[u_1] \qquad \forall t \in (0, 1).
\]
Hence, $J: \cE\to \R$ is strictly concave.
\end{proof}

\begin{theorem}[An example of a non-concave and non-convex reaction functional with the nonlinear kinetics]
\label{t:MM}
Let $\Omega=\{x\in\R^3:1<|x|<2\}$,
 $f(s)=s/(1+s)$ $(s > -1),$
$\sigma=1/50$, $\kappa=1/5$, and $G=e^4-1$. Define
\[
g(x)=
\frac{G\sinh\bigl(\kappa(2-|x|)\bigr)}
{|x|\sinh\kappa} \quad \mbox{and} \quad 
u(x)=\frac{1+g(x)}{50} \quad \forall x \in \overline{\Omega}, 
\quad \mbox{and} \quad M_{\rm e}=\int_\Omega u\,dx>0.
\]
Then $u \in \cE$ and $\rho_u:=\cT_g(u) =g$. Moreover,  there exist 
$v,w\in\cA_u$ such that
\[
J''[u][v]>0
\qquad\mbox{and}\qquad
J''[u][w]<0.
\]
Consequently, $J$ is neither concave nor convex on $\cE$.
\end{theorem}

\begin{proof}

We choose $u$ and $g$ so that \eqref{BVP4rho} has the
explicit solution $\rho_u=g$. To do so, we require 
$uf (\rho_u) = \sigma g$, i.e., $u f(g) = \sigma g$. 
Since $f(s) = s/(s+1)$, this leads to 
$u = \sigma (1 + g)$, as we assumed in the theorem. 
Now \eqref{BVP4rho} becomes 
\begin{equation}
\label{Onlyg}
\left\{
\begin{aligned}
& -\Delta g+2\sigma g =0 \quad\mbox{in }\Omega,\\
&g=G\quad\mbox{on }\{|x|=1\} \quad \mbox{and} \quad 
g=0\quad\mbox{on }\{|x|=2\}.
\end{aligned}
\right.
\end{equation}
We verify that the assumed expression of $g$ is the radial (i.e., radially symmetric, function of $|x|$) solution of
this problem. Thus, $g$ satisfies \eqref{BVP4rho},
and uniqueness in Proposition~\ref{p:FirstPaper}
gives $\rho_u=g$ in ${\Omega}$.

We  first construct $w \in \cA_u$ such that $J''[u][w] < 0.$
Let $w(x)=x_1$ for $x=(x_1,x_2,x_3)\in\overline\Omega$.
Clearly $w \in C(\overline{\Omega})$ and by symmetry the integral of $w$ over $\Omega$ vanishes. 
Since  $u\ge1/50$, it is verified that $w \in \cA_u.$

Let $\psi_u = \cL_u(uf'(\rho_u))= \cL_u(uf'(g))$ (cf.\ \eqref{cLu} and \eqref{BVPpsi}) and set $q=1-\psi_u$. 
By Proposition~\ref{p:psiu}, $0 < q \le 1$
on $\overline{\Omega}.$
Since $g$ and $u$ are smooth and radial, 
elliptic regularity
(cf.\ Theorem~8.13 in \cite{GilbargTrudinger98}) and uniqueness
imply that $q$ is smooth and radial.
Write $\rho = \rho_u = g$ in $\Omega.$
Let $\hat{\psi}_{u,w}
=\cL_u\bigl(-wf(\rho_u)\bigr)
=-\cL_u\bigl(wf(\rho)\bigr).$ By \eqref{cLu}, 
\[
\left\{ 
\begin{aligned}
-\Delta\hat{\psi}_{u,w}
+(\sigma+uf'(\rho))\hat{\psi}_{u,w}
&=-x_1f(\rho)&&\quad \mbox{in }\Omega,\\
\hat{\psi}_{u,w}&=0&&\quad \mbox{on }\partial\Omega.
\end{aligned}
\right.
\]
Define $\tilde{\psi}(x)
=-\hat{\psi}_{u,w}(-x_1,x_2,x_3)$.
Since all of $\Omega$, $u$, and $\rho = g$ are radially symmetric, $\tilde{\psi}$ satisfies the same boundary value problem as $\hat{\psi}_{u,w}$.
By uniqueness, $\tilde{\psi}=\hat{\psi}_{u,w}$,
so $\hat{\psi}_{u,w}$ is odd in $x_1$. Hence $\hat{\psi}_{u, w} = 0$ on $\Omega \cap \{ x_1 = 0\}$
and furthermore, $\hat{\psi}_{u,w}=0$ on
$\partial(\Omega\cap\{x_1>0\})$.

Since $\rho > 0$ in $\Omega$ (cf.\ the definition of $g$), we have $-x_1f(\rho)<0$ in
$\Omega\cap\{x_1>0\}$.
The weak maximum principle then implies $\hat{\psi}_{u,w}\le0$
there. Since $\hat{\psi}_{u,w}\not\equiv0$ there, the strong maximum principle gives
$\hat{\psi}_{u,w}<0$ in $\Omega\cap\{x_1>0\}$ (cf.\ Theorems~8.1 and 8.19 in \cite{GilbargTrudinger98}).
By oddness, $w\hat{\psi}_{u,w}<0$ whenever $x_1\ne0$.
Since $q>0$, $u > 0$, $f'(\rho)>0$, and $f''(\rho)<0$ in $\Omega,$ 
the formula for the second variation in \eqref{J''uv} implies 
\begin{equation}
\label{J''uw}
J''[u][w]
=\int_\Omega q\left[
2wf'(\rho)\hat{\psi}_{u,w}
+uf''(\rho)\hat{\psi}_{u,w}^2
\right]dx<0.
\end{equation}

We now construct $v \in \cA_u$ and prove $J''[u][v]>0$ in three steps. 

\textit{Step 1. A reformulation of $J''[u][v]$.} We first recall from \eqref{J''uv} the formula
of the second variation
\begin{equation}
\label{RecallJpp}
J''[u][v]
=
2\int_\Omega qvf'(\rho)\hat{\psi}_{u,v}\,dx
+\int_\Omega quf''(\rho)\hat{\psi}_{u,v}^2\,dx
\qquad\forall v\in\cA_u.
\end{equation}
Since $f''(\rho)<0$, the second integral is non-positive. We therefore need to choose $v$ so that the first integral
is positive and sufficiently large. The equation for $\hat{\psi}_{u,v}=\cL_u ( - v f(\rho))$ is
\begin{equation}
\label{minusvfrho}
-\Delta\hat{\psi}_{u,v}
+(\sigma+uf'(\rho))\hat{\psi}_{u,v}
=-vf(\rho).
\end{equation}
To use this equation in the first integral, we seek a function
$h$ such that $\displaystyle f'(\rho)\hat{\psi}_{u,v}=-f(\rho)h$, so that the 
first integral in \eqref{RecallJpp} becomes
\[
-2\int_\Omega qh\,vf(\rho)\,dx,
\] 
in which $vf(\rho)$
can be replaced using \eqref{minusvfrho} and $\hat{\psi}_{u,v} = -hf(\rho)/f'(\rho).$  Since $f(s) = s/(s+1)$, 
we have 
$f(\rho)/f'(\rho)=\rho(1+\rho)$. So, the required relation is
\begin{equation}
\label{hat}
\hat{\psi}_{u,v}=-\rho(1+\rho)h.
\end{equation}
Thus, we seek a smooth function $h$, vanishing on
$\partial\Omega$, and a direction $v\in\cA_u$ satisfying 
\eqref{minusvfrho},  which by \eqref{hat} is 
\begin{equation}
\label{vfrho}
vf(\rho)
=-\Delta\bigl(\rho(1+\rho)h\bigr)
+(\sigma+uf'(\rho))\rho(1+\rho)h.
\end{equation}

We now derive an identity for $J''[u][v]$, and then 
choose $h$ and verify that
\eqref{vfrho} defines an admissible direction $v \in \cA_u.$
Since $f'(\rho)=(1+\rho)^{-2}$ and
$f''(\rho)=-2(1+\rho)^{-3}$, 
we get by 
\eqref{RecallJpp}, \eqref{hat}, and \eqref{vfrho} that 
\begin{align*}
J''[u][v]
&=-2\int_\Omega qvf(\rho)h\,dx
-2\int_\Omega qu\frac{\rho^2}{1+\rho}h^2\,dx\\
&=2\int_\Omega qh\Delta\bigl(\rho(1+\rho)h\bigr)\,dx
-2\sigma\int_\Omega q\rho(1+\rho)h^2\,dx
\\
&\qquad -2\int_\Omega qu\left[\rho(1+\rho)f'(\rho)
+\frac{\rho^2}{1+\rho}\right]h^2\,dx\\
&=2\int_\Omega qh\Delta\bigl(\rho(1+\rho)h\bigr)\,dx
-2\int_\Omega q\bigl(\sigma\rho(1+\rho)+u\rho\bigr)h^2\,dx.
\end{align*}
We next rewrite the first integral 
so that we can use the equations for $\rho$ and $q$. Since $h=0$ on $\partial\Omega$
by our choice of $h$, integration by parts gives
\begin{align*}
&2\int_\Omega qh\Delta\bigl(\rho(1+\rho)h\bigr)\,dx\\
&\qquad 
=-2\int_\Omega \nabla(qh)\cdot\nabla\bigl(\rho(1+\rho)h\bigr)\,dx\\
&\qquad =-2 \int_\Omega \left[ h \nabla q + q \nabla h\right] \cdot \left[ h\nabla ( \rho (1+\rho)) + 
\rho (1 + \rho) \nabla h\right] dx
\\
&\qquad =-2\int_\Omega q\rho(1+\rho)|\nabla h|^2\,dx
-2\int_\Omega \left[q\nabla\bigl(\rho(1+\rho)\bigr)
+\rho(1+\rho)\nabla q\right]\cdot h\nabla h\,dx\\
&\qquad \quad-2\int_\Omega h^2\nabla q\cdot\nabla\bigl(\rho(1+\rho)\bigr)\,dx.
\end{align*}
To use the equations for $\rho$ and $q$, we convert
the mixed derivative terms into Laplacians. In the second integral, the expression in brackets is
$\nabla\bigl(q\rho(1+\rho)\bigr)$ and 
$2h\nabla h=\nabla(h^2)$. Thus, integrating by parts
once more, we obtain from the previous equation that 
\begin{align}
\label{DDD}
&2\int_\Omega qh\Delta\bigl(\rho(1+\rho)h\bigr)\,dx\nonumber \\
&=-2\int_\Omega q\rho(1+\rho)|\nabla h|^2\,dx
+\int_\Omega\left[\Delta\bigl(q\rho(1+\rho)\bigr)
-2\nabla q\cdot\nabla\bigl(\rho(1+\rho)\bigr)\right]h^2\,dx\nonumber \\
&=-2\int_\Omega q\rho(1+\rho)|\nabla h|^2\,dx
+\int_\Omega\left[q\Delta\bigl(\rho(1+\rho)\bigr)
+\rho(1+\rho)\Delta q\right]h^2\,dx,
\end{align}
where the last equality follows by expanding
$\Delta\bigl(q\rho(1+\rho)\bigr)$.
Thus,
\begin{align}
\label{EEE}
J''[u][v]
&=-2\int_\Omega q\rho(1+\rho)|\nabla h|^2\,dx\nonumber \\
&\quad+\int_\Omega\Bigl[
q\Delta\bigl(\rho(1+\rho)\bigr)
+\rho(1+\rho)\Delta q
-2q\bigl(\sigma\rho(1+\rho)+u\rho\bigr)
\Bigr]h^2\,dx.
\end{align}
We recall the equations for $\rho = \rho_u$ (cf.\ \eqref{BVP4rho})
and $q = 1 - \psi_u$ (cf.\ \eqref{BVPpsi}), and expand $\Delta \bigl( \rho (1 + \rho) \bigr)$,
\begin{align*}
&\Delta\rho=\sigma\rho+uf(\rho), \\
& \Delta q=(\sigma+uf'(\rho))q-\sigma, \\
&\Delta\bigl(\rho(1+\rho)\bigr)
=(1+2\rho)\Delta\rho+2|\nabla\rho|^2.
\end{align*} 
With these, 
we can simplify the expression in the brackets in the second integral of \eqref{EEE}
\begin{align}
\label{FFF}
&q\Delta\bigl(\rho(1+\rho)\bigr)
+\rho(1+\rho)\Delta q
-2q\bigl(\sigma\rho(1+\rho)+u\rho\bigr)\nonumber\\
&\qquad =2q|\nabla\rho|^2
+uq\left[(1+2\rho)f(\rho)+\rho(1+\rho)f'(\rho)-2\rho\right]
\nonumber\\
&\qquad \quad +\sigma q\left[\rho(1+2\rho)-\rho(1+\rho)\right]
-\sigma\rho(1+\rho)\nonumber\\
&\qquad =2q|\nabla\rho|^2+\sigma q\rho^2-\sigma\rho(1+\rho),
\end{align}
where the last equality follows from $f(\rho) = \rho/(1+\rho)$ and 
$f'(\rho) = (1+\rho)^{-2}$, and hence
\[
(1+2\rho)f(\rho)+\rho(1+\rho)f'(\rho)-2 \rho = 0,
\]
eliminating $u$ from the equation. Substituting \eqref{FFF} into the expression \eqref{EEE} for $J''[u][v]$ gives
\begin{equation}
\label{J''h}
J''[u][v]
=2\int_\Omega q
|\nabla\rho|^2h^2 \, dx- 2\int_\Omega q\rho  (1+\rho)|\nabla h|^2 \, dx -\sigma\int_\Omega\rho(1+\rho-q\rho)h^2\,dx.
\end{equation}

\textit{Step 2. Construction of $v$ in $\cA_u$.} 
We seek $h$ for which the first integral in \eqref{J''h}
exceeds the other two combined. To satisfy the zero-integral condition on $v$, we take
$h$ to be $x_1/|x|$ times a function of $\rho$.
Since $\rho$ and the coefficients of $\cL_u^{-1}$ are
radial, the relation \eqref{vfrho} then makes $v$ odd
in $x_1$. 

To choose that function, set $t=\log(1+\rho)$.
Since $\nabla\rho=(1+\rho)\nabla t$, including the factor
$(1+\rho)^{-1}$ in $h$ cancels the factor $(1+\rho)^2$
in $|\nabla\rho|^2h^2$, which appears in the first
integral of \eqref{J''h}. We therefore try
$h=(x_1/|x|)e^{-t}P(t)$, where $P$ is a polynomial
to be chosen. Using a polynomial allows us to impose
the required boundary behavior through its zeros and
to carry out explicit estimates in Step~3.
We next choose $P$ so that $h$ vanishes on the boundary
and $v$ extends smoothly to $\overline\Omega$.
Since $t=4$ on the inner sphere and $t=0$ on the
outer sphere, we require $P(0)=P(4)=0$.
We also impose $P'(0)=0$ to ensure that $v$ remains
smooth where $\rho=0$. Indeed, this makes
$\rho(1+\rho)h$ contain a factor $\rho^3$, so applying
$\cL_u^{-1}$ leaves a factor $\rho$ that cancels
the denominator $f(\rho)$.
The simplest choice is $t^2(4-t)$. 

To that end, we define
\begin{align*}
&P(t)=t^2(4-t) \qquad \forall t\in\R,\\
&h(x)=\frac{x_1}{|x|}
\frac{P(\log(1+\rho(x)))}{1+\rho(x)}
\qquad \forall x\in\overline{\Omega},\\
&v=\frac{\cL_u^{-1}\bigl(\rho(1+\rho)h\bigr)}{f(\rho)}
\qquad\mbox{in }\Omega.
\end{align*}
Since $h=0$ on $\partial\Omega$,
$\rho(1+\rho)h\in H_0^1(\Omega)\cap H^2(\Omega)$.
Thus, $\cL_u^{-1}\bigl(\rho(1+\rho)h\bigr)$ is well defined;
cf.\ \eqref{cLu}.
Define $\alpha(s)=\log(1+s)/s$ for $s>-1$, $s\ne0$,
and $\alpha(0)=1$. Then
$\alpha\in C^\infty((-1,\infty))$.
We have
\begin{equation}
\label{ax}
\rho(1+\rho)h=x_1\rho^3a
\qquad \mbox{with} \qquad 
a(x)=\frac{\alpha(\rho(x))^2}{|x|}
\bigl[4-\log(1+\rho(x))\bigr].
\end{equation}
In particular, $a\in C^\infty(\overline{\Omega})$.
It then follows from the definition of $v$ and that of $\cL_u$ (cf.\ \eqref{cLu}), $f(\rho) = \rho/(1+\rho)$, 
and \eqref{ax} that 
\begin{align*}
v (x) & = \frac{1}{f(\rho(x))} \left[ \cL_u^{-1} (\rho (1+\rho )h)\right] (x)\\
& = \frac{1}{f(\rho (x))} \left[ - \Delta \left( x_1 \rho^3(x) a(x) \right) 
+ (\sigma + u(x)f'(\rho (x)) ) x_1 \rho^3(x) a(x) \right]\\
& = -x_1  b(x) \qquad \forall x \in \Omega,
\end{align*}
where 
\[
\begin{aligned}
b=(1+\rho)\Bigl[
&\rho^2\Delta a
+6a|\nabla\rho|^2
+3a\rho\Delta\rho
+6\rho\nabla\rho\cdot\nabla a\\
&-(\sigma+uf'(\rho))a\rho^2
+2\rho^2\frac{\partial_{x_1}a}{x_1}
+6a\rho\frac{\partial_{x_1}\rho}{x_1}
\Bigr].
\end{aligned}
\]
For $r=|x|$, the radial symmetry implies 
$\partial_{x_1}a/x_1=a'(r)/r
$ and $ \partial_{x_1}\rho/x_1=\rho'(r)/r.$
At $x_1=0$, these quotients are understood through
the smooth extensions given by the right-hand sides.
Thus, $b$ is smooth and radial on $\overline{\Omega}$.
Consequently, $v=-x_1b$ extends smoothly to
$\overline{\Omega}$ and is odd in $x_1$, so
by symmetry, its integral vanishes on $\Omega$. Since $v$ is bounded and $u\ge1/50$, we have
$v\in\cA_u$.

\textit{Step 3. Positive second variation.} We now compare the positive and negative contributions
in \eqref{J''h}.
Since $\rho$ and $q$ are radial, we use spherical coordinates $r = |x|.$
We define 
\begin{align*}
&m(r)=-r^2\rho'(r) \qquad  \forall r \in [1, 2], \\
&t(r)=\log(1+\rho(r)) \qquad \forall r \in [1, 2], \\
&A=\int_0^4 e^{-t}P(t)^2\,dt,\\
&
Q=\int_0^4(e^{-t}-e^{-2t})
\bigl(P(t)^2-P'(t)^2\bigr)\,dt.
\end{align*}
We make the following three claims.

\textit{Claim 1.} We have
\begin{align}
\frac{299}{300}\le q\le1
&\qquad\mbox{on }\overline{\Omega},
\label{qm}\\
\frac{149}{150}\,2G\le m(r)\le\frac{151}{150}\,2G
&\qquad\forall r\in[1,2].
\label{mr}
\end{align}

\textit{Claim 2.} We have
\begin{equation}
\label{AQP}
0<A<17,\qquad Q>\frac52,\qquad
0\le P(t)^2<90\quad\forall t\in[0,4].
\end{equation}

\textit{Claim 3.} With
$B=\frac{299}{300}\frac{149}{150}A
-\frac{151}{150}(A-Q)$,
the second variation satisfies
\begin{equation}
\label{J''uvB}
J''[u][v]\ge\frac{4\pi}{3}
\left[
4GB-4\int_1^2P(t(r))^2\,dr
-\sigma\int_1^2r^2P(t(r))^2\,dr
\right].
\end{equation}

Assuming these claims, the definition of $B$ and
\eqref{AQP} give
\[
B =\frac{151}{150}Q-\frac{749}{45000}A
\ge Q-\frac{A}{50} >\frac52-\frac{17}{50}>2.
\]
Since $0\le\rho(r)\le G=e^4-1$, we have
$0\le t(r)\le4$ on $[1,2]$.
Thus, \eqref{AQP} gives $P(t(r))^2<90$.
Using this and $\sigma=1/50$, the integrals in \eqref{J''uvB} satisfy
\begin{align*}
4\int_1^2P(t(r))^2\,dr
+\sigma\int_1^2r^2P(t(r))^2\,dr
&<360+\frac{90}{50}\cdot\frac73<365.
\end{align*}
Finally, with $B>2$ and $G=e^4-1>53$ in \eqref{J''uvB}, we conclude that
\begin{equation}
\label{J''uvG}
J''[u][v]>\frac{4\pi}{3}(8G-365)>0.
\end{equation}

We now prove these three claims.

\textit{Proof of Claim 1.}
To estimate $q$, as $uf'(\rho)=1/[50(1+\rho)]\le1/50$,
we compare $\psi_u$ with the solution of
\[
-\Delta p=1/50\quad\mbox{in }\Omega,
\qquad p=0\quad\mbox{on }\partial\Omega.
\]
Solving this gives
$p(x)=(7-|x|^2-6/|x|)/300$.
Moreover, for $1 \le |x| \le 2$, the maximum of $p$ is attained at $|x| = 3^{1/3}.$
Thus, 
$0\le p\le(7-3^{5/3})/300<1/300$.
Since $u=(1+\rho)/50$,
\[
\cL_u^{-1}p
=\frac1{50}+(\sigma+uf'(\rho))p
\ge\frac1{50}
\ge\frac1{50(1+\rho)}
=uf'(\rho)
=\cL_u^{-1}\psi_u.
\]
The comparison principle gives $0\le\psi_u\le p$.
Thus, since $p < 1/300,$
\[
 299/300< 1 - p \le q=1-\psi_u\le1,
\]
proving
\eqref{qm}.

Since $\rho$ is radial and $\kappa^2=2\sigma=1/25$, \eqref{Onlyg} gives $\Delta\rho=\kappa^2\rho$,
which is $\rho''(r) + (2/r) \rho'(r) = \kappa^2 \rho.$ Hence, $m'(r)=-\kappa^2r^2\rho(r)\le0$.
Differentiating the explicit formula for $\rho= g$ gives
\begin{equation}
\label{mmm}
\frac{2G\kappa}{\sinh\kappa}
=m(2)\le m(r)\le m(1)
=G(1+\kappa\coth\kappa).
\end{equation}
By Taylor's expansion, 
\[
\frac{\sinh a}{a} = \sum_{n=0}^\infty \frac{a^{2n}}{(2n+1)!} < \sum_{n=0}^\infty \left(\frac{a^2}{6}\right)^n
= \frac{1}{1 - a^2/6} \qquad \mbox{if } 0 < a < \sqrt{6}. 
\]
Since $\kappa = 1/5< \sqrt{6}$, we then have 
\begin{equation}
\label{kappa1}
\frac{\kappa}{\sinh\kappa} \ge 1-\frac{\kappa^2}{6} = \frac{149}{150}. 
\end{equation}
Define 
\[
\zeta (a) = \left( 1  + \frac{a^2}{3} \right) \sinh a - a \cosh a \qquad \forall a \ge0. 
\]
Then, we have $\zeta(0) = 0$ and also $\zeta'(a) > 0$ for $a>0$. Hence, 
\begin{equation}
\label{kappa2}
\kappa \coth \kappa < 1+\frac{\kappa^2}{3} = \frac{76}{75}
\end{equation}
Now \eqref{mr} follows from \eqref{mmm}, \eqref{kappa1} and \eqref{kappa2}.

\textit{Proof of Claim 2.}
Direct integration gives
\[
A=144-6992e^{-4},
\qquad
Q=\frac{721}{8}-4728e^{-4}-\frac{953}{8}e^{-8}.
\]
Since $54<e^4<55$,
\[
0<A<144-\frac{6992}{55}<17,
\qquad
Q>\frac{721}{8}-\frac{4728}{54}
-\frac{953}{8\cdot54^2}>\frac52.
\]
Also, $P'(t)=t(8-3t)$, so
$0\le P(t)\le P(8/3)=256/27$ on $[0,4]$.
Consequently, $P(t)^2\le(256/27)^2<90$.
This proves \eqref{AQP}.

\textit{Proof of Claim 3.}
Write $x=r\omega$, where $\omega = (\omega_1, \omega_2, \omega_3)\in\mathbb S^2$.
By the definition of $h$,
$h(r\omega)=\omega_1e^{-t(r)}P(t(r))$.
Also,
\begin{equation}
\label{tprime}
t'(r)=\frac{\rho'(r)}{1+\rho(r)}
=-\frac{m(r)}{r^2(1+\rho(r))}<0.
\end{equation}
Thus $t$ decreases from $4$ to $0$ as $r$ increases
from $1$ to $2$. We denote its inverse by $r(t)$.
By symmetry with respect to $\omega_i $ $ (i = 1, 2, 3)$, 
\[
\int_{\mathbb S^2}\omega_1^2\,dS=\frac{4\pi}{3}.
\]
Since $\rho$ and $q$ are radial, integration in spherical coordinates and applying the change of variables
from $r$ to $t$ by $dr=-r(t)^2e^t\,dt/m(r(t))$, together with the definition 
$m(r) = - r^2 \rho'(r)$, \eqref{qm}, \eqref{mr}, and the definition of $A$, imply 
that the first integral in $J''[u][v]$ in \eqref{J''h} is
\begin{align*}
\int_\Omega q|\nabla\rho|^2h^2\,dx
&=\frac{4\pi}{3}\int_1^2
q(r)r^2\rho'(r)^2e^{-2t(r)}P(t(r))^2\,dr\\
&=\frac{4\pi}{3}\int_0^4
m(r(t))q(r(t))e^{-t}P(t)^2\,dt\\
&\ge\frac{4\pi}{3}\,2G
\frac{149}{150}\frac{299}{300}A.
\end{align*}

We next estimate the second integral of $J''[u][v]$ in \eqref{J''h}. 
To estimate $|\nabla h|^2$, we  first separate 
the radial
and tangential contributions of $\nabla h$. Since $r=|x|$ and $\omega_1=x_1/r$, we have
$
\nabla r=\omega,
$ and $ \nabla\omega_1=(e_1-\omega_1\omega)/{r}, $
where $e_1=(1,0,0)$.
Differentiating $h(r\omega)=\omega_1e^{-t(r)}P(t(r))$ gives
\begin{align*}
\nabla h(r\omega)
&=\omega_1e^{-t(r)}
\bigl[P'(t(r))-P(t(r))\bigr]t'(r)\omega +\frac{e^{-t(r)}P(t(r))}{r}
(e_1-\omega_1\omega).
\end{align*}
Since 
$\omega\cdot(e_1-\omega_1\omega)=0$
and $|e_1-\omega_1\omega|^2
=1-\omega_1^2, $
we have
\begin{align}
\label{Gradh2}
|\nabla h(r\omega)|^2
&=\omega_1^2e^{-2t(r)}
\bigl[ P'(t(r))-P(t(r))\bigr]^2t'(r)^2+\frac{1-\omega_1^2}{r^2}
e^{-2t(r)}P(t(r))^2.
\end{align}
Since
\begin{align*}
&\int_{\mathbb S^2}(1-\omega_1^2)\,dS=\frac{8\pi}{3}, \\
&\rho(1+\rho) = f(\rho) (1+\rho)^2 = f(\rho)e^{2 t},
\end{align*}
we obtain by \eqref{Gradh2}, \eqref{tprime} for a formula of $t'(r)$, the change of variable  
$dr=-r(t)^2e^t\,dt/m(r(t))$, and the definition $t(r) = \log (1 + \rho(r))$ that
\begin{align*}
&\int_\Omega\rho(1+\rho)|\nabla h|^2\,dx\\
&\qquad =\frac{4\pi}{3}\int_1^2
r^2f(\rho(r))t'(r)^2
\bigl(P'(t(r))-P(t(r))\bigr)^2\,dr
+\frac{8\pi}{3}\int_1^2 f(\rho(r))P(t(r))^2\,dr\\
&\qquad =\frac{4\pi}{3}\int_0^4
m(r(t))(e^{-t}-e^{-2t})
\bigl(P'(t)-P(t)\bigr)^2\,dt
+\frac{8\pi}{3}\int_1^2 f(\rho(r))P(t(r))^2\,dr.
\end{align*}
Since $P(0)=P(4)=0$, integration by parts gives
\begin{align*}
&\int_0^4(e^{-t}-e^{-2t})
\bigl(P'(t)-P(t)\bigr)^2\,dt\\
&\qquad = \int_0^4(e^{-t}-e^{-2t})
\bigl(P'(t)^2+P(t)^2\bigr)\,dt
-\int_0^4(e^{-t}-e^{-2t})(P(t)^2)'\,dt\\
&\qquad =\int_0^4(e^{-t}-e^{-2t})P'(t)^2\,dt
+\int_0^4e^{-2t}P(t)^2\,dt\\
&\qquad =A-Q.
\end{align*}
It follows from $q\le1$, $f(\rho)\le1$, and
\eqref{mr} that we obtain the estimate for the second integral 
in $J''[u][v]$ in \eqref{J''h}: 
\begin{align*}
\int_\Omega q\rho(1+\rho)|\nabla h|^2\,dx
&\le\int_\Omega\rho(1+\rho)|\nabla h|^2\,dx\\
&\le\frac{4\pi}{3}
\left[
2G\frac{151}{150}(A-Q)
+2\int_1^2P(t(r))^2\,dr
\right].
\end{align*}

Finally, we note that 
\[
\frac{\rho(1+\rho-q\rho)}{(1+\rho)^2}
=f(\rho)(1-qf(\rho))\in[0,1].
\]
Thus, for the last integral in $J''[u][v]$ in \eqref{J''h}, we have 
\begin{align*}
\sigma\int_\Omega\rho(1+\rho-q\rho)h^2\,dx
&=\frac{4\pi\sigma}{3}\int_1^2
r^2f(\rho(r))
\bigl(1-q(r)f(\rho(r))\bigr)P(t(r))^2\,dr\\
&\le\frac{4\pi\sigma}{3}\int_1^2
r^2P(t(r))^2\,dr.
\end{align*}
Substituting these three estimates into
\eqref{J''h} proves \eqref{J''uvB}.

Finally, by Proposition~\ref{p:variations},
$s\mapsto J[u+sz]$ is $C^2$ near $s=0$ for $z=v,w$.
Hence \eqref{J''uw} and \eqref{J''uvG}
exclude convexity and concavity, respectively.
\end{proof}

\section{Non-Existence of Critical Points}
\label{s:CriticalPoints}

We recall the definition of a critical point
(cf.\ Definition~5.1 in \cite{KhanLiFirst}).

\begin{definition}
\label{d:CriticalPoint}
We say that $u\in\cE$ is a critical point of the reaction
functional $J:\cE\to\R$ if $J'[u][v]=0$ for all $v\in\cA_u.$
\end{definition}

\begin{theorem}[Non-existence of critical points]
\label{t:CriticalPoint}
For both the linear and nonlinear cases, the reaction functional
$J:\cE\to\R$ admits no critical points. More precisely,
let $u\in\cE$ and define
\begin{equation}
\label{vu}
v_u=u\left[f(\rho_u)(1-\psi_u)
-\frac1{M_{\rm e}}\int_\Omega u f(\rho_u)(1-\psi_u)\,dx\right]
\qquad \mbox{in } \Omega. 
\end{equation}
Then $v_u\in\cA_u$ and $J'[u][v_u]>0$.
\end{theorem}

\begin{proof}
Let $u\in\cE$. We first prove that $v_u \in \cA_u. $ 
Recall from Proposition~\ref{p:FirstPaper} that $\rho_u = \cT_g(u) \in H^1_g(\Omega)\cap H^2(\Omega)$ and 
from Proposition~\ref{p:psiu} that $\psi_u=\cL_u\bigl(u f'(\rho_u)\bigr) \in H^1_0(\Omega)\cap H^2(\Omega).$
Denote
\[
c_u=\frac1{M_{\rm e}}\int_\Omega u (1-\psi_u) f(\rho_u)\,dx.
\]
Note that $u \ge 0$ a.e.\ in $\Omega$ and 
$| \{ u > 0 \}| > 0$, since $u \in \cE.$ 
Moreover, $0 < \rho_u < g_{\rm max} $ and $ 0 \le \psi_u < 1$ in $\Omega$ by Proposition~\ref{p:FirstPaper} and 
Proposition~\ref{p:psiu}, respectively. Hence, 
\[
0\le f(\rho_u) (1 - \psi_u) \le f(g_{\rm max}) \qquad \mbox{in } \Omega, 
\]
and 
$0<c_u\le f(g_{\rm max})$.  Consequently, it follows 
from the definition of $v_u$ (cf.\ \eqref{vu})
and the fact that $u \ge 0$ a.e.\ in $\Omega$ that
\begin{align*}
|v_u| & \le u \left| \left[f(\rho_u)(1-\psi_u)
-\frac1{M_{\rm e}}\int_\Omega u f(\rho_u)(1-\psi_u)\,dx\right]  \right| 
\\
& \le u  \left[  f(g_{\rm max}) + f(g_{\rm max} ) \frac{1}{M_{\rm e}} \int_\Omega u \, dx \right]  \\
& \le 2 f(g_{\rm max}) u \qquad \mbox{a.e. in } \Omega. 
\end{align*}
Moreover, by the definition of $v_u$ (cf.\ \eqref{vu}), 
\begin{equation*}
    \int_\Omega v_u\,dx=
\int_\Omega u(1-\psi_u) f(\rho_u)\,dx-c_uM_{\rm e}=0.
\end{equation*}
Therefore, setting $t_v = ( 2 f(g_{\rm max}))^{-1} > 0$, we have $u + t v_u \in \cE$ 
for all $t \in (-t_v, t_v)$ and thus $v_u\in\cA_u$.

Now, we show that $J'[u][v_u] > 0.$ By \eqref{J'uv} and \eqref{vu}, and noting that
the integral of $u$ over $\Omega$ is $M_{\rm e}, $ we obtain 
\begin{align}
\label{J'vu}
J'[u][v_u]
&=\int_\Omega u (1-\psi_u) f(\rho_u)\left[ (1-\psi_u) f(\rho_u)-c_u\right] \,dx
\nonumber \\
&=\int_\Omega u \left[ (1-\psi_u) f(\rho_u)-c_u\right]^2\,dx
+ c_u \int_\Omega u \left[ f(\rho_u) (1 - \psi_u) - c_u \right] dx\nonumber \\
& = \int_\Omega u \left[ (1-\psi_u) f(\rho_u)-c_u\right]^2\,dx \nonumber \\
&\ge 0.
\end{align}
We prove that the equality is impossible. Suppose otherwise. Then, by the third
equality in \eqref{J'vu}, 
\begin{equation}
\label{cu}
(1-\psi_u) f(\rho_u)=c_u\qquad\mbox{a.e.\ in }\{u>0\}.
\end{equation}
Denote $q = 1 - \psi_u \in H^2(\Omega)$. Note by Proposition~\ref{p:psiu} that $q=1$ on $\partial \Omega$
and $q > 0$ in $\overline{\Omega}$. Setting
\begin{equation} 
\label{K}
K=\frac{1+g_{\rm max}}{\min_{\overline\Omega}q} >0, 
\end{equation}
we claim that $q$ satisfies
\begin{equation}
\label{Deltaqc}
\Delta q -K|\nabla q |^2\le0\qquad\mbox{a.e.\ in }\Omega. 
\end{equation}

Suppose \eqref{Deltaqc} is true. We show there will be a contradiction. This in turn will show that $J'[u][v_u] > 0$. 
Set $w=-e^{-Kq}\in H^2(\Omega)$. By
\eqref{Deltaqc} and the fact that $q=1$ on $\partial\Omega$, we have
\begin{align*}
&\Delta w=Ke^{-Kq}(\Delta q-K|\nabla q|^2)\le0
\qquad\mbox{a.e.\ in }\Omega,\\
&w=-e^{-K}\qquad\mbox{on }\partial\Omega.
\end{align*}
The weak minimum principle
(cf.\ Theorem~8.1 in \cite{GilbargTrudinger98}) implies
$w\ge-e^{-K}$ on $\overline{\Omega}$. On the other hand,
$q \le1$ on $\overline{\Omega}$ implies $w=-e^{-Kq }\le-e^{-K}$ on $\overline{\Omega}$. It then follows that
$w=-e^{-K}$ and $q =1$ throughout $\Omega$. Substitution in \eqref{BVPpsi} now gives
$uf'(\rho_u)=0$ a.e.\ in $\Omega$. 
For both of the reaction kinetics, linear and nonlinear (cf.\ \eqref{f}), 
\[
f'(\rho_u)\ge \frac{1}{(1+g_{\rm max})^{2}}>0 \qquad \mbox{in }\Omega.
\]
Hence $u=0$ a.e.\ in $\Omega$. This contradicts the fact that $u \in \cE$ which implies 
\[
\int_\Omega u\,dx=M_{\rm e}>0.
\]

We now prove our claim, i.e., \eqref{Deltaqc}, using \eqref{cu}. 
Denote for simplicity $\rho = \rho_u$ and recall that $ q = 1- \psi_u.$
Since  $\rho, q \in H^2(\Omega)$ by Proposition~\ref{p:FirstPaper} and Proposition~\ref{p:psiu}
and all $f$, $f'$, and $f''$ are bounded on $[0,g_{\rm max}]$,  
we have $qf(\rho)\in H^2(\Omega)$.
Note that the weak derivatives of any $H^1(\Omega)$-function vanish on its level-set
(cf.\ Lemma~7.7 in \cite{GilbargTrudinger98}). By \eqref{cu} and two applications of this result, 
$ \nabla ( qf (\rho)) = 0$ a.e.\ in $\{u>0\}$, i.e., 
\begin{equation}
    \label{nablarho}
\nabla\rho=-\frac{f(\rho)}{qf'(\rho)}\nabla q \qquad \mbox{a.e. in } \{ u > 0 \},
\end{equation}
and 
\begin{equation}
    \label{Deltaqfrho}
\Delta(qf(\rho))=\nabla \cdot \nabla (q f(\rho)) = 0 \qquad \mbox{a.e.\ in } \{u>0\}. 
\end{equation}
Note by Proposition~\ref{p:FirstPaper} and Proposition~\ref{p:psiu} that
$q > 0$ and $f(\rho) > 0$ in $\Omega.$
The equations for $\rho$ and $q$ (cf.\ \eqref{BVP4rho} and \eqref{BVPpsi}) lead to
\begin{align}
\label{Laplacerho}
&\Delta\rho=\sigma\rho+uf(\rho) \qquad \mbox{in } \Omega,\\
\label{Laplaceq}
&\Delta q=\sigma(q-1)+uf'(\rho)q \qquad \mbox{in } \Omega.
\end{align}
These identities imply that the following
holds true a.e.\ in $\{u>0\}$: 
\begin{align}
\label{qfrho}
0 & = \Delta (q f(\rho))  \qquad \mbox{[by \eqref{Deltaqfrho}]}
\nonumber \\
&=f(\rho)\Delta q+qf'(\rho)\Delta\rho
+qf''(\rho)|\nabla\rho|^2
+2f'(\rho)\nabla q\cdot\nabla\rho 
\nonumber \\
&=f(\rho)\Delta q+\sigma q\rho f'(\rho)
+uqf'(\rho)f(\rho)
+\frac{f(\rho)^2f''(\rho)}{qf'(\rho)^2}|\nabla q|^2
-\frac{2f(\rho)}{q}|\nabla q|^2
\nonumber \\
&\qquad \qquad \qquad \quad \ \mbox{[by \eqref{Laplacerho} for $\Delta \rho$ and 
\eqref{nablarho} for $\nabla \rho$]}
\nonumber \\
&=2f(\rho)\Delta q -\sigma f(\rho )(q -1)
+\sigma q \rho f'(\rho)
+\frac{f(\rho)}{q}
\left[\frac{f(\rho)f''(\rho)}{f'(\rho)^2}-2\right]|\nabla q|^2.
\\
&\qquad \qquad \qquad \quad \ \mbox{[by \eqref{Laplaceq} for $\Delta q$]}
\nonumber 
\end{align}
For the linear case, $f(\rho)=\rho$, $f'(\rho)=1$, and
$f''(\rho)=0$. For the nonlinear case,
\[
\frac{\rho f'(\rho)}{f(\rho)}=\frac1{1+\rho}
\qquad \mbox{and} \qquad 
\frac{f(\rho)f''(\rho)}{f'(\rho)^2}=-2\rho.
\]
Hence, it follows from \eqref{qfrho} that
\begin{equation*}
\Delta q=
\begin{cases}
-\dfrac\sigma2+\dfrac{|\nabla q|^2}{q}
&\mbox{for the linear case},\\
\dfrac\sigma2\left(q-1-\dfrac {q}{1+\rho}\right)
+\dfrac{1+\rho}{q}|\nabla q|^2
&\mbox{for the nonlinear case},
\end{cases}
\qquad \mbox{a.e. in } \{ u > 0 \}.
\end{equation*}
Since $0<q\le1$ and
$0<\rho\le g_{\rm max}$ in $\Omega$ and $\sigma > 0$, 
with our definition of the constant $K$ (cf.\ \eqref{K}), both cases imply 
\begin{equation}
    \label{K1}
\Delta q \le\frac{1+g_{\rm max}}{\min_{\overline\Omega}q }
|\nabla q |^2= K | \nabla q|^2\qquad\mbox{a.e.\ in } \{u>0\}.
\end{equation}
It follows from \eqref{Laplaceq} that
\begin{equation}
\label{K2}
\Delta q=\sigma(q-1)\le0 \qquad \mbox{a.e. on } \{ u = 0 \}.
\end{equation}
Now \eqref{Deltaqc} follows from \eqref{K1} and \eqref{K2}. 
\end{proof}

\section{Characterization of maximizing sequences}
\label{s:Characterization}

We first present some basic properties of maximizing sequences of the functional 
$J: \cE\to \R$, and then show sufficient and necessary conditions for such sequences and discuss the 
related special and simple sufficient conditions. 

Let $\rho_0=\cT_g(0)$ be the substrate concentration in the
absence of enzymes, i.e., $\rho_0 \in H_g^1(\Omega)\cap H^2(\Omega)$ is the unique solution to the 
boundary-value problem
\begin{equation}
   \label{BVPrho0}
    \left\{
    \begin{aligned}
         - \Delta \rho_0 + \sigma \rho_0 &= 0 & & \quad \mbox{in } \Omega, \\
         \rho_0 &= g  & & \quad \mbox{on } \partial \Omega. 
    \end{aligned}
    \right.
\end{equation}
We have $0 \le \rho_0 \le g_{\rm max}$ on $\overline{\Omega}$ by the maximum principle. 
Recall that 
\begin{equation*} 
\cL_0=(-\Delta+\sigma)^{-1}: L^2(\Omega) \longrightarrow H_0^1(\Omega) \cap H^2(\Omega) 
\end{equation*}
is the operator $\cL_u$ with $u  = 0$ in $\Omega;$ cf.\ \eqref{cLu}. 
If $v \in L^2(\Omega)$, then $w := \cL_0 v \in H_0^1(\Omega)\cap H^2(\Omega)$ is the unique solution
to the boundary-value problem
\begin{equation}
  \label{BVPw} 
    \left\{
    \begin{aligned}
         - \Delta w + \sigma w &= v & & \quad \mbox{in } \Omega, \\
         w &= 0 & & \quad \mbox{on } \partial \Omega. 
    \end{aligned}
    \right.
\end{equation}
It is clear that $\cL_0$ is linear and self-adjoint. Moreover, by the
 maximum principle, it is positive, i.e.,  $\cL_0 v \ge 0$ if $v\ge 0$ in $\Omega. $

\begin{lemma}
    \label{l:L0}
We have 
        \begin{equation}
\label{vL0v}
\int_\Omega v\cL_0 v\,dx
=\int_\Omega\left(
|\nabla\cL_0 v|^2+\sigma|\cL_0v|^2
\right)dx
\qquad\forall v\in L^2(\Omega).
\end{equation}
Moreover, if we denote by $ \| v \|_\sigma $ the square root of the common value from the above identity, 
then $\| \cdot \|_\sigma $ is a norm on $L^2(\Omega)$, and it is equivalent to the $H^{-1}(\Omega)$-norm. 
\end{lemma}

\begin{proof}
 Let $v \in L^2(\Omega)$ and $w = \cL_0v \in H_0^1(\Omega)\cap H^2(\Omega).$ 
Let $\hat{v} \in L^2(\Omega)$. Setting the test function to be 
$\hat{w} = \cL_0 \hat{v}$ in the weak formulation for \eqref{BVPw}, we obtain
\begin{equation}
\label{innersigma}
\langle v,\hat v\rangle_\sigma := 
\int_\Omega v \cL_0 \hat{v} \, dx =  \int_\Omega 
\left( \nabla \cL_0 v \cdot
\nabla \cL_0 \hat{v}   + \sigma \cL_0 v \, \cL_0 \hat{v} \right) dx. 
\end{equation}
This defines an inner
product on $L^2(\Omega)$. Taking $\hat v=v$ gives
\eqref{vL0v}, and the corresponding norm is
$\|\cdot\|_\sigma$.

Let $v \in L^2(\Omega)$ and $\xi \in H_0^1(\Omega)$. From the weak formulation for \eqref{BVPw}, we get
\begin{align*}
\int_\Omega v \xi \, dx &= \int_\Omega \left( \nabla \cL_0 v \cdot \nabla \xi + 
\sigma \cL_0 v \, \xi \right) dx \\
& \le  \| \nabla \cL_0 v \|_{L^2(\Omega)} \| \nabla \xi \|_{L^2(\Omega)}
+ \sigma \| \cL_0 v \|_{L^2(\Omega)} \| \xi \|_{L^2(\Omega)} \\
& \le (1+\sqrt{\sigma} )\| v \|_{\sigma} \| \xi \|_{H^1(\Omega)}. 
\end{align*}
Hence, $\| v \|_{H^{-1}(\Omega)} \le (1+\sqrt{\sigma}) \| v \|_{\sigma}.$ 
For the converse, we may assume $v\ne0$,
since the case $v=0$ is immediate. Setting $\xi = \cL_0 v \in H^1_0(\Omega)$, we get by \eqref{vL0v} that
\[
\| v \|_{H^{-1}(\Omega)} \ge \frac{1}{\| \cL_0 v \|_{H^1(\Omega)} } \int_\Omega v \cL_0 v \, dx 
\ge \sqrt{\min\{ 1, \sigma\} }\| v \|_\sigma, 
\]
completing the proof.
\end{proof}

The following estimates will be used to characterize maximizing
sequences. Part~{\rm (2)} gives sufficient conditions, while
Part~{\rm (3)} provides the converse estimates by restricting $u$
to the set $\{\rho_u\ge \theta g_{\rm max}\}$ for any $0 < \theta <1$. 
Note that $\rho_u = \cT_g(u)\in H_g^1(\Omega)\cap H^2(\Omega)$ is defined in Proposition~\ref{p:FirstPaper}.
We denote by $\1_A$ the indicator function  of a set $A.$

\begin{lemma}
\label{l:JuvL0}
\begin{compactenum}
\item[\rm (1)] We have for any $u \in L^2_+(\Omega)$ that
        \begin{equation}
\label{rho0rhou}
0\le\rho_0-\rho_u
=\cL_0 \bigl(uf(\rho_u)\bigr)
\le f(g_{\rm max})\cL_0 u
\qquad \mbox{in } \Omega. 
\end{equation}
\item[{\rm (2)}]
If $u\in \cE$, $v\in L^2(\Omega)$, and $0\le v\le u$ a.e.\ in
$\Omega$, then \begin{align}
\label{JuvL0}
0 & < M_{\rm e}f(g_{\rm max})-J[u] \nonumber \\
& \le \int_\Omega
\bigl[ f(g_{\rm max})-f(\rho_0)\bigr] u\,dx
+f(g_{\rm max})\left(
\|u-v\|_{L^1(\Omega)} + \| v\|_\sigma^2 \right).
\end{align}

\item[{\rm (3)}]
If $u \in \cE$, $0 < \theta < 1$, and  $v_\theta =u\1_{\{\rho_u\ge \theta g_{\rm max}\}}$, then
\begin{align}
\label{uvL1}
&\|u-v_\theta \|_{L^1(\Omega)}
\le
\frac{M_{\rm e}f(g_{\rm max})-J[u]}
{f(g_{\rm max})-f(\theta g_{\rm max})},\\
\label{vL0vJ}
&
\| v_\theta\|_\sigma^2 
\le
\frac{(1+g_{\rm max})^2}{f(\theta g_{\rm max})}
\bigl( M_{\rm e}f(g_{\rm max})-J[u]\bigr) .
\end{align}
\end{compactenum}
\end{lemma}

\begin{proof}
(1) The inequality $\rho_0 \ge \rho_u$ in $\Omega$ follows from Proposition~\ref{p:FirstPaper}. 
Subtracting the equation for $\rho_u$ from that for $\rho_0$ leads to the equality in 
\eqref{rho0rhou} and applying then the comparison principle implies the last inequality in \eqref{rho0rhou}. 

(2) By \eqref{supJ}, $M_{\rm e} f(g_{\rm max})> J[u]. $
We now estimate the difference between the reaction at
$\rho_0$ and that at $\rho_v$. Since $0<f'(s)\le1$ on
$[0,g_{\rm max}]$, it follows from \eqref{rho0rhou} (with $v$ replacing $u$) 
and the definition of $\| v\|_\sigma $ in Lemma~\ref{l:L0} that
\begin{equation}
\label{Jvrho0}
\begin{aligned}
0\le\int_\Omega vf(\rho_0)\,dx-J[v]
&\le\int_\Omega v(\rho_0-\rho_v)\,dx\le f(g_{\rm max})\int_\Omega v\cL_0 v\,dx
= f(g_{\rm max}) \| v \|_\sigma^2. 
\end{aligned}
\end{equation}
Since $u \ge v$ a.e.\ in $\Omega$, 
Proposition~\ref{p:FirstPaper} (3) implies that $J[u]\ge J[v]$. Thus, since
$u \in \cE$, $\rho_0 \le g_{\rm max}$ in $\Omega$,  and $J[u] <  M_{\rm e} f(g_{\rm max})$,
we obtain by \eqref{Jvrho0} that
\begin{align*}
0 &< M_{\rm e}f(g_{\rm max})-J[u]\\
& \le M_{\rm e}f(g_{\rm max})-J[v]\\
&=f(g_{\rm max})\int_{\Omega} (u-v) \, dx
+f(g_{\rm max})\int_\Omega v\,dx-J[v]\\
& =f(g_{\rm max})\|u-v\|_{L^1(\Omega)}
+\int_\Omega v\bigl[ f(g_{\rm max})-f(\rho_0)\bigr]\,dx
+\int_\Omega vf(\rho_0)\,dx-J[v]\\
&\le  \int_\Omega \bigl[ f(g_{\rm max})-f(\rho_0)\bigr] u\,dx
+f(g_{\rm max})\left( \|u-v\|_{L^1(\Omega)}
+ \| v \|_\sigma^2\right), 
\end{align*}
proving \eqref{JuvL0}.

(3)
Since $f(\rho_u)\le f(\theta g_{\rm max})$ on $ \{ \rho_u < \theta g_{\rm max} \}$, we have
\begin{align*}
\bigl[ f(g_{\rm max})-f(\theta g_{\rm max})\bigr]
\|u-v_\theta \|_{L^1(\Omega)}
& = \int_{\{ \rho_u < \theta g_{\rm max} \} } u \bigl[ f(g_{\rm max})-f(\theta g_{\rm max})\bigr]\, dx
\\
&\le 
\int_{\{\rho_u<\theta g_{\rm max}\}}
u\bigl[ f(g_{\rm max})-f(\rho_u)\bigr]\,dx\\
& \le
\int_\Omega u\bigl[f(g_{\rm max})-f(\rho_u)\bigr]\,dx
\\
& =M_{\rm e}f(g_{\rm max})-J[u].
\end{align*}
This proves \eqref{uvL1}.

Note that
$ uf (\rho_u) \ge v_\theta f(\rho_u) \ge v_\theta f(\theta g_{\rm max})$ a.e. in $\Omega$. 
Positivity and linearity of $\cL_0: L^2(\Omega)\longrightarrow H_0^1(\Omega)\cap H^2(\Omega)$ 
together with \eqref{rho0rhou} then lead to 
\[
\rho_0-\rho_u
=\cL_0 \bigl(uf(\rho_u)\bigr) 
\ge \cL_0\bigl(v_\theta f(\theta g_{\rm max})\bigr)
= f(\theta g_{\rm max})\cL_0 v_\theta \qquad \mbox{in } \Omega. 
\]
Multiplying by $v_\theta $ and integrating and noting that $0 \le \rho_0 \le g_{\rm max}$ on $\overline{\Omega}$,  we obtain
\[
f(\theta g_{\rm max})\int_\Omega v_\theta \cL_0 v_\theta \,dx
\le\int_\Omega v_\theta (\rho_0-\rho_u)\,dx
\le\int_\Omega u(g_{\rm max}-\rho_u)\,dx.
\]
For either reaction law in \eqref{f},
$f'(s)\ge(1+g_{\rm max})^{-2}$ on $[0,g_{\rm max}]$.
Consequently,
\begin{equation}
\label{fgmaxfs}
g_{\rm max}-s
\le(1+g_{\rm max})^2
\bigl[ f(g_{\rm max})-f(s)\bigr] 
\qquad\forall s\in[0,g_{\rm max}].
\end{equation}
It follows from the above two inequalities together with the definition 
of $\| v_\theta \|_\sigma^2$ in Lemma~\ref{l:L0} that
\begin{align*}
f(\theta g_{\rm max}) \| v_\theta \|_{\sigma}^2 
&= 
f(\theta g_{\rm max})\int_\Omega v_\theta \cL_0  v_\theta \,dx\\
&\le(1+g_{\rm max})^2
\int_\Omega u\bigl[ f(g_{\rm max})-f(\rho_u)\bigr]\,dx\\
&=(1+g_{\rm max})^2
\bigl( M_{\rm e}f(g_{\rm max})-J[u]\bigr),
\end{align*}
which leads to \eqref{vL0vJ}.
\end{proof}

Before we characterize the maximizing sequences of the functional $J: \cE\to \R$, we present in 
the next proposition some properties of such sequences. We denote
\begin{equation}
    \label{Gammamax}
\Gamma_{\rm max}=\{y\in\partial\Omega:g(y)=g_{\rm max}\}.
\end{equation}

\begin{proposition}
\label{p:rhokuk}
\begin{compactenum}
\item[{\rm (1)}]
Let $\{u_k\}\subset\cE$ be a maximizing sequence of $J: \cE \to \R$. Then
$\|u_k \|_{L^2(\Omega)} \to \infty$ as $k \to \infty.$
\item[{\rm (2)}]
Let $K\subset\overline\Omega\setminus\Gamma_{\rm max}$ be
nonempty and compact. We have 
\begin{equation}
\label{uK}
0 \le \int_{K}u\,dx
\le\frac{M_{\rm e}f(g_{\rm max})-J[u]}
{\min_K(f(g_{\rm max})-f(\rho_0))}\qquad \forall u \in \cE.
\end{equation}
Consequently, if $\{ u_k \} $ is a maximizing sequence of $J: \cE \to\R$, then
\begin{equation}
\label{0KO}
0 \le \int_{K}u_k\,dx \to 0 \qquad \mbox{as } k \to \infty.
\end{equation}
\item[{\rm (3)}]
We have
\begin{equation}
\label{rhou0}
\|\rho_u-\rho_0\|_{H^1(\Omega)}^2
\le\frac{f(g_{\rm max})(1+g_{\rm max})^2}{\min\{1,\sigma\}}
\bigl( M_{\rm e}f(g_{\rm max})-J[u]\bigr)
\qquad\forall u\in\cE.
\end{equation}
Consequently, if $\{u_k\}\subset\cE$ is a maximizing
sequence, then
$
\rho_{u_k}\to\rho_0$ in $ H^1(\Omega).
$
\end{compactenum}
\end{proposition}

\begin{proof}
(1) Suppose that $\|u_k\|_{L^2(\Omega)}$ does not tend to
$+\infty$. Then a subsequence, not relabeled, is bounded
in $L^2(\Omega)$. Passing to a further subsequence, we have
$u_k\wc u$ in $L^2(\Omega)$. Hence, $u \in \cE$ and $J[u_k]\to J[u]$ by Proposition~\ref{p:FirstPaper}.
Since $\{ u_k \}$ is a maximizing sequence, $J[u] = \max_{v \in \cE}J[v]$. 
This contradicts \eqref{supJ}, the non-attainment of the supremum of $J: \cE \to \R.$

(2) Let $u \in \cE.$ Note that 
$
\min_K \left[ f(g_{\rm max})-f(\rho_0) \right] >0
$
and $\rho_0 \ge \rho_u$ in $\Omega$ by Proposition~\ref{p:FirstPaper}.
Hence, 
\[
\begin{aligned}
\min_K \left[ f(g_{\rm max})-f(\rho_0) \right] 
\int_{K}u\,dx
&\le\int_\Omega
\bigl[ f(g_{\rm max})-f(\rho_0)\bigr] u\,dx 
\le M_{\rm e}f(g_{\rm max})-J[u].
\end{aligned}
\]
This proves \eqref{uK}. If $\{ u_k \}$ is a maximizing sequence, then
$J[u_k] \to M_{\rm e} f(g_{\rm max})$. This and \eqref{uK} imply \eqref{0KO}. 

(3)
Let $u\in\cE$ and set
$w=\rho_0-\rho_u\in H_0^1(\Omega)$.
By \eqref{rho0rhou} and \eqref{fgmaxfs},
\[
0\le w\le g_{\rm max}-\rho_u
\le(1+g_{\rm max})^2
\bigl[ f(g_{\rm max})-f(\rho_u)\bigr] \qquad \mbox{in }\Omega. 
\]
Testing the equation for $w$ against $w$, we obtain
\[
\begin{aligned}
\min\{1,\sigma\}\|w\|_{H^1(\Omega)}^2
&\le\int_\Omega\bigl(|\nabla w|^2+\sigma w^2\bigr)\,dx\\
&=\int_\Omega u f(\rho_u)w\,dx\\
&\le f(g_{\rm max})(1+g_{\rm max})^2
\int_\Omega u\bigl[ f(g_{\rm max})-f(\rho_u)\bigr] \,dx\\
&=f(g_{\rm max})(1+g_{\rm max})^2
\bigl(M_{\rm e}f(g_{\rm max})-J[u]\bigr).
\end{aligned}
\]
This proves \eqref{rhou0}, and the stated convergence as $M_{\rm e} f(g_{\rm max})
= \sup_{\cE} J.$
\end{proof}

\begin{theorem}[Characterization of maximizing sequences]
\label{t:ukJ}
Suppose $\{u_k\} \subset \cE$. Then $\{u_k\}$ is a maximizing sequence of the reaction
functional $J: \cE \to \R$ if
and only if the following hold true: 
\begin{compactenum}
\item[\rm (1)]
\begin{equation}
\label{ukrho0}
\int_\Omega
\bigl[ f(g_{\rm max})-f(\rho_0)\bigr] u_k\,dx\to0 \qquad \mbox{as } k\to \infty; 
\end{equation}
\item[\rm (2)]
there exist $v_k\in L^2(\Omega)$ $(k = 1, 2, \dots)$ such that $0\le v_k\le u_k$ a.e.\ in $\Omega$ for each $k \ge 1$ 
and that 
\begin{equation}
\label{ukvk}
\|u_k-v_k\|_{L^1(\Omega)}\to0
\qquad \mbox{and} \qquad 
\| v_k \|_{H^{-1}(\Omega)} \to 0 \qquad \mbox{as } k \to \infty.
\end{equation}
\end{compactenum}
\end{theorem}

\begin{proof}
Let $u\in\cE$. By Proposition~\ref{p:FirstPaper}, $\rho_u\le\rho_0\le g_{\rm max}$ in $\Omega.$
\begin{equation}
\label{ukJ}
0\le\int_\Omega u 
\bigl[ f(g_{\rm max})-f(\rho_0)\bigr] \, dx
\le M_{\rm e}f(g_{\rm max})-J[u].
\end{equation}
If $\{u_k\}$ is a maximizing sequence, then \eqref{ukJ} implies
\eqref{ukrho0}. Choosing
$
v_k=u_k\1_{\{\rho_{u_k}\ge g_{\rm max}/2\}},
$
Lemma~\ref{l:JuvL0} with $\theta=1/2$, together
with the norm equivalence in Lemma~\ref{l:L0},
implies \eqref{ukvk}.

Conversely, if \eqref{ukrho0} and \eqref{ukvk} hold, then Lemma~\ref{l:L0} gives $\|v_k\|_\sigma\to0$. Applying \eqref{JuvL0} to $u_k$ and $v_k$ yields  
$J[u_k]\to M_{\rm e}f(g_{\rm max})$.
\end{proof}

Setting $v_k = u_k$ $(k \ge 1)$ leads to some sufficient conditions for $\{ u_k \}
\subset \cE$ to be a maximizing sequence. The following theorem indicates that 
the maximizing sequence constructed in our previous work \cite{KhanLiFirst} satisfies such conditions and in general
such conditions are not necessary conditions:

\begin{theorem}[Sufficient conditions for maximizing sequences]
\label{t:sufficient}
\begin{compactenum}
\item[{\rm (1)}] If $\{ u_k \}\subset \cE$ satisfies
\begin{equation}
\label{sufficient}
\int_\Omega
\bigl[f(g_{\rm max})-f(\rho_0)\bigr]u_k\,dx\to0
\qquad\mbox{and}\qquad
\| u_k \|_{H^{-1}(\Omega)}\to 0
\qquad \mbox{as } k\to \infty,
\end{equation}
then it is a maximizing sequence. Moreover, there exists a sequence $\{ u_k \} \subset \cE$ that satisfies
these sufficient conditions and hence is a maximizing sequence. 
\item[{\rm (2)}]
The sufficient conditions \eqref{sufficient} are not necessary conditions. More precisely, there 
exists a maximizing sequence $\{u_k\}\subset\cE$ such that $\|u_k\|_{H^{-1}(\Omega)}\to \infty$. 
In particular, the auxiliary sequence $\{v_k\}$ in
Theorem~\ref{t:ukJ} cannot in general be replaced by $\{u_k\}$.
\end{compactenum}
\end{theorem}

\begin{remark}
    If $\{ u_k \} \subset \cE$ is a maximizing sequence of the reaction functional $J: \cE \to \R$, 
    then Proposition~\ref{p:rhokuk} gives $\| u_k \|_{L^2(\Omega) } \to \infty$. 
    By Lemma~\ref{l:L0} and Theorem~\ref{t:sufficient}, we may have $\|u_k \|_{H^{-1}(\Omega)}
    \to 0$ or $\| u_k \|_{H^{-1}(\Omega)} \to \infty.$
\end{remark}

\begin{proof}[Proof of Theorem~\ref{t:sufficient}]
(1) 
The sufficiency of the conditions \eqref{sufficient} follows from Theorem~\ref{t:ukJ} with 
$v_k = u_k$ $(k\ge 1).$ We now show that the maximizing sequence $\{ u_k \} \subset \cE$
constructed in our previous work \cite{KhanLiFirst}
satisfies these sufficient conditions, i.e., we can choose $v_k = u_k $ $(k\ge 1)$ in Theorem~\ref{t:ukJ}.

Fix $y_0\in\partial\Omega$ with $g(y_0)=g_{\rm max}$. Let $\ve > 0$ be small and denote 
$\Gamma_\ve = \partial \Omega \cap B_\ve(y_0).$
Define the normal coordinate map $F(z,t)=z-t\nu(z)$ for $z \in \Gamma_\ve$ and $t \in (0, \ve^3)$.
Define finally
\[
S_\ve=F\bigl(\Gamma_\ve\times(0,\ve^3)\bigr)\qquad \mbox{and} \qquad 
u_\ve=\frac{M_{\rm e}}{|S_\ve|}\1_{S_\ve}.
\]
Note that $u_\ve$ is the function constructed in \cite{KhanLiFirst}. Moreover, it is shown
in \cite{KhanLiFirst} that, for any given 
 $\ve_k\downarrow0$, the sequence
$u_k=u_{\ve_k}$ $(k\ge 1)$  is a maximizing sequence of the reaction functional $J: \cE\to \R.$
Here, we show $\{ u_k \}$ satisfies the sufficient conditions \eqref{sufficient}. 

For sufficiently small $\ve>0$, the Jacobian of $F$ is
bounded above and below by positive constants independent
of $\ve$. Moreover, there exist constants $c>0$ and $C > 0$ such that  
\[
c\ve^5\le|S_\ve|\le C\ve^5 \qquad
\mbox{and} \qquad S_\ve\subset B_{\ve+\ve^3}(y_0).
\]
Let $\rho_0 = \cT_g(0)$ be given by \eqref{BVPrho0}. 
Since $\rho_0$ is continuous on $\overline\Omega$ and
$\rho_0(y_0)=g_{\rm max}$,
\[
0\le\int_\Omega
\bigl[f(g_{\rm max})-f(\rho_0)\bigr]u_\ve\,dx
\le M_{\rm e}\sup_{S_\ve}
\bigl[f(g_{\rm max})-f(\rho_0)\bigr]\to0.
\]
Thus, the first part of the conditions \eqref{sufficient} holds true for $u_k = u_{\ve_k}$ $(k \ge 1).$

Let $\varphi\in C_{\rm c}^\infty(\Omega)$. Since $\varphi=0$
on $\partial\Omega$, it follows that $\varphi(F(z,0)) = 0$ for any $z \in \Gamma_\ve$. 
Applying the Cauchy--Schwarz inequality, we obtain 
for $z\in\Gamma_\ve$ and $0<t<\ve^3$ that 
\begin{align*}
|\varphi(F(z,t))|^2&= \left|\int_0^t  \partial_s \bigl( \varphi (F(z, s))\bigr)\, ds \right|^2 \\
&
=\left|\int_0^t\nabla\varphi(F(z,s))\cdot\nu(z)\,ds\right|^2\\
&\le\ve^3\int_0^{\ve^3}|\nabla\varphi(F(z,s))|^2\,ds.
\end{align*}
Integrating, using the change of variable $x = F(z, t)$ and the bounds on the Jacobian, we get for some constant $C > 0$ that
\begin{align*}
\int_{S_\ve}|\varphi(x)|^2\,dx
&\le C\int_{\Gamma_\ve}\int_0^{\ve^3}
|\varphi(F(z,t))|^2\,dt\,dS(z)\\
&\le C\ve^6\int_{\Gamma_\ve}\int_0^{\ve^3}
|\nabla\varphi(F(z,s))|^2\,ds\,dS(z)\\
&\le C\ve^6\|\nabla\varphi\|_{L^2(\Omega)}^2.
\end{align*}
By density, this estimate holds for every
$\varphi\in H_0^1(\Omega)$. Consequently, 
the definition of $u_\ve$, the Cauchy--Schwarz inequality, and 
the lower bound $|S_\ve | \ge c\ve^5$ imply
\[
\left|\int_\Omega u_\ve\varphi\,dx\right| = \frac{M_{\rm e}}{|S_\ve|}  \left|\int_{S_\ve}  \varphi \, dx \right|
\le\frac{M_{\rm e}}{|S_\ve|^{1/2}}
\|\varphi\|_{L^2(S_\ve)}
\le CM_{\rm e}\ve^{1/2}\|\varphi\|_{H^1(\Omega)},
\]
and hence
$\|u_\ve\|_{H^{-1}(\Omega)}\le CM_{\rm e}\ve^{1/2}\to0$. With $\ve = \ve_k$ and 
$u_k = u_{\ve_k}$, this implies the second
part of \eqref{sufficient}.

(2) Let $\{\alpha_k\}\subset\cE$ 
be a maximizing sequence. 
(For instance, we can take the maximizing sequence constructed above.)
We construct a new maximizing sequence
$\{ u_k \} \subset \cE$ such that $\| u_k \|_{H^{-1}(\Omega)} \to \infty.$
Fix  $x_0\in\Omega$. Choose $\delta_k\downarrow0$ and
$r_k\downarrow0$ such that
\[
\delta_k<M_{\rm e},\qquad
B(x_0,2r_k)\subset\Omega,\qquad
\frac{\delta_k^2}{r_k}\to\infty.
\]
Choose nonnegative $\beta_k\in C^\infty_{\rm c}(B(x_0,r_k))$
with 
\[
\int_\Omega \beta_k\,dx=\delta_k \qquad \forall k \ge 1.
\]
Finally, set
\[
u_k=\left(1-\frac{\delta_k}{M_{\rm e}}\right)\alpha _k+\beta_k \qquad \forall k \ge 1.
\]
Then $u_k\in\cE$ $(k \ge 1)$ and 
\[
\|u_k-\alpha_k\|_{L^1(\Omega)}\le2\delta_k\to 0 
\qquad \mbox{as } k\to \infty.
\]
Proposition~\ref{p:Juv} therefore implies
$J[u_k]-J[\alpha_k]\to0$, so $\{u_k\}$ is also a maximizing
sequence.

Now, fix  $\phi\in C^\infty_{\rm c}(B(0,2))$ with
$\phi=1$ on $B(0,1)$ and set
$\phi_k(x)=\phi((x-x_0)/r_k)$ for each $k \ge 1$. Then
$\phi_k\in C^\infty_{\rm c}(B(x_0,2r_k))$ and
$\phi_k=1$ on $B(x_0,r_k)$.
Since $B(x_0,2r_k)\subset\Omega$, the change of variables
$x=x_0+r_ky$ gives
\begin{equation}
\label{Crk}
\int_\Omega\bigl(|\nabla\phi_k|^2+\sigma\phi_k^2\bigr)\,dx
=r_k\int_{\R^3}|\nabla\phi|^2\,dy
+\sigma r_k^3\int_{\R^3}\phi^2\,dy
\le Cr_k,
\end{equation}
where
\[
C=\int_{\R^3}|\nabla\phi|^2\,dy
+\sigma r_1^2\int_{\R^3}\phi^2\,dy > 0
\]
is independent of $k$, since $0 < r_k\le r_1$ for all $k \ge 1.$

Testing the equation for $\cL_0\beta _k$ against $\phi_k$, we obtain by
 \eqref{vL0v}, the Cauchy--Schwarz inequality, and \eqref{Crk} that
\[
\delta_k
=\int_\Omega \beta_k\phi_k\,dx
=\int_\Omega
\left[
\nabla(\cL_0 \beta_k)\cdot\nabla\phi_k
+\sigma(\cL_0 \beta_k)\phi_k
\right]dx
\le
\left(\int_\Omega \beta_k\cL_0  \beta_k\,dx\right)^{1/2} (Cr_k)^{1/2}.
\]
Since $u_k\ge \beta_k\ge0$, the positivity of $\cL_0$ implies
$\cL_0 u_k\ge\cL_0 \beta_k\ge0$. Hence,
\[
\| u_k \|_\sigma^2 = \int_\Omega u_k\cL_0 u_k\,dx
\ge\int_\Omega \beta_k\cL_0  \beta_k\,dx
\ge\frac{\delta_k^2}{Cr_k}\to\infty \qquad \mbox{as } k \to \infty.
\]
By Lemma~\ref{l:L0}, $\| u_k \|_{H^{-1}(\Omega)} \to \infty.$ 
\end{proof}

\section{Limiting Radon measures}
\label{s:LimitMeasures}

Let $\mathcal M(\overline{\Omega})$ be the Banach space of signed Radon measures on $\overline{\Omega}$ \cite{folland}.
If $\mu \in \mathcal M(\overline{\Omega})$, the norm is the total variation $\| \mu \|=|\mu|(\overline{\Omega}).$
If $u \in L^1(\Omega)$, then it is the density of a Radon measure $\mu \in \mathcal M(\overline{\Omega})$. 
This is denoted by $u dx = d \mu$.  
If $\nu \in \mathcal M(\overline{\Omega}) $ and $v_k \in L^1(\Omega)$ $(k \ge 1)$ with each $v_k$ the density
of $\nu_k \in \mathcal M(\overline{\Omega})$, then the weak-$\ast$ convergence
$\nu_k \stackrel{*}{\rightharpoonup} \nu$ will also be denoted by 
$v_k\,dx\stackrel{*}{\rightharpoonup}d \nu$ on $\overline\Omega$.
This is defined by
\begin{equation}
\label{varphivknu}
\int_{\overline{\Omega}} v_k\varphi\,dx
\to
\int_{\overline\Omega}\varphi\,d\nu
\qquad\forall\varphi\in C(\overline\Omega).
\end{equation}
In this case, we also say that the sequence $v_kdx $ converges weakly-$\ast$  to $d \nu$ or to $\nu$
on $\overline{\Omega}.$
We denote 
\begin{align}
\label{MathcalMmax}
&
\mathcal M_{\rm max}(\overline{\Omega}) = \{ \mu \in \mathcal M(\overline{\Omega}): \mu \ge 0, 
\ \| \mu \| = M_{\rm e}, \mbox{ and } \mbox{supp}\,(\mu) 
\subseteq \Gamma_{\rm max} \},
\end{align}
where $\Gamma_{\rm max}$ is defined in \eqref{Gammamax}. 

The following is our main result of this section; it characterizes the limiting
measures of maximizing sequences of the functional $J: \cE \to \R$: 

\begin{theorem}[Limiting Radon measures]
\label{t:RadonMeasures}
\begin{compactenum}
\item[{\rm (1)}]
Let $\{u_k\}\subset\cE$ be a maximizing sequence of $J: \cE\to \R$. Then
$\{u_k\,dx\}$ has a weak-$\ast$  convergent subsequence. Moreover, every subsequential
weak-$\ast$  limiting measure of $\{ u_k \, dx \}$ is a measure in $\mathcal M_{\rm max}(\overline{\Omega}).$ 
\item[{\rm (2)}]
Let $\mu \in \mathcal M_{\rm max}(\overline{\Omega})$ and $0 \le \ell \le M_{\rm e} f(g_{\rm max}) = 
\sup_{\cE} J.$ There exist  $ u_k\in \cE \cap C_{\rm c}^\infty(\Omega)$ $(k = 1, 2, \dots)$ such that 
$u_k\,dx\stackrel{*}{\rightharpoonup}d\mu $ and 
$J[u_k] \to \ell$.
\end{compactenum}
\end{theorem}

\begin{remark}
   Part (2) of the theorem indicates that a limiting measure $\mu \in \mathcal M_{\rm max}(\overline{\Omega})$ 
   does not determine the limiting value of the total reaction flux (i.e., the functional) $J$. 
\end{remark}

The following are immediate consequences of Theorem~\ref{t:RadonMeasures}: 

\begin{corollary}
\label{c:Dirac}
Suppose $\Gamma_{\rm max}=\{y_0\}$ and
$\{u_k\}\subset\cE$ is a maximizing sequence of $J: \cE \to \R$. Then
$
u_k\,dx\stackrel{*}{\rightharpoonup}M_{\rm e}\delta_{y_0}$ on $\overline\Omega.$
\end{corollary}

\begin{proof}
The sequence $\{u_k\,dx\}$ is relatively weak-$\ast$  compact.
By Theorem~\ref{t:RadonMeasures}, every subsequential limit is a
nonnegative measure of mass $M_{\rm e}$ supported on
$\{y_0\}$, and hence equals $M_{\rm e}\delta_{y_0}$.
Thus the original sequence converges weakly-$\ast$ to $M_{\rm e} \delta_{y_0}.$
\end{proof}

\begin{corollary}
    \label{c:mu} 
    A nonnegative measure $\mu\in \mathcal M (\overline{\Omega})$ is the
weak-$\ast$ limit of $\{u_k\,dx\}$ for some maximizing sequence
$\{u_k\}\subset\cE$ of $J: \cE \to \R$ if and only if $\mu \in \mathcal M_{\rm max}(\overline{\Omega}).$
\end{corollary}

\begin{proof}
    The ``if" part follows from Part (2) and the ``only if" part follows from Part (1) of Theorem~\ref{t:RadonMeasures}.
\end{proof}

The proof of Theorem~\ref{t:RadonMeasures} uses some lemmas which provide detailed constructions
of needed sequences of enzyme concentrations. 

\begin{proof}[Proof of Theorem~\ref{t:RadonMeasures}]
(1) The measures $u_k\,dx$ are nonnegative and have mass
$M_{\rm e}$ on the compact space $\overline\Omega$.
Identifying $\mathcal M(\overline\Omega)$ with
$C(\overline\Omega)^*$, the Banach--Alaoglu theorem and
the separability of $C(\overline\Omega)$ give a weak-$\ast$
convergent subsequence of $\{u_k\,dx\}$
(cf.\ Corollary~7.18 and Theorem~5.18 in \cite{folland}).

Suppose a subsequence $\{ u_{k_j}dx \} $ converges weakly-$\ast$ to some $\mu \in \mathcal M(\overline{\Omega})$.
By \eqref{varphivknu} with $u_{k_j}$ replacing $v_k$, we get $\mu \ge 0.$ 
Moreover, choosing  $\varphi = 1$ and sending $j\to \infty$, we obtain $\| \mu \| = M_{\rm e}.$
Applying \eqref{ukJ} to $u = u_{k_j}$ for each $j$ and sending $j\to \infty$, we obtain
\begin{equation}
\label{fgfr}
\int_{\overline\Omega}
\bigl[ f(g_{\rm max})-f(\rho_0)\bigr] \,d\mu=0.
\end{equation}
Note that $\alpha := f(g_{\rm max}) - f(\rho_0) \in C(\overline{\Omega})$, 
$\alpha \ge 0$ on $\overline{\Omega}$, and $\alpha > 0$ on $\overline{\Omega}\setminus
\Gamma_{\rm max}$. Therefore, $\mu (\overline{\Omega}\setminus \Gamma_{\rm max}) = 0$. 
Since $\overline{\Omega}\setminus \Gamma_{\rm max}$ is open in $\overline{\Omega}$, we have
$\mbox{supp}\,(\mu) \subset \Gamma_{\rm max}.$

(2) If $\mu \in \mathcal M_{\rm max}(\overline\Omega)$, then $\mbox{supp}\,(\mu) \subseteq \Gamma_{\rm max}$, 
$\mu(\Gamma_{\rm max}) = \mu(\overline{\Omega})
= M_{\rm e}$,  and 
\[
\int_{\partial\Omega}f(g)\,d\mu = f(g_{\rm max}) \mu (\Gamma_{\rm max})
=M_{\rm e}f(g_{\rm max}).
\]
Part (2) of Lemma~\ref{l:mu} below gives a sequence
$\{u_k\}\subset\cE$ such that
$u_k\,dx\stackrel{*}{\rightharpoonup}d\mu$ and
$J[u_k]\to\ell$ for the prescribed
$\ell\in[0,M_{\rm e}f(g_{\rm max})]$.
\end{proof}

We now prove two lemmas. The first one will be used to prove the second one, 
which was used in proving Theorem~\ref{t:RadonMeasures} above.

\begin{lemma}
\label{l:construction}
There exist $r_0>0$ and $C>0$, depending only on $\Omega$,
such that, given $y\in\partial\Omega$, $0<h<r<r_0$, 
and $b>0$, there exists a nonnegative $v\in C^\infty_{\rm c}(\Omega)$
satisfying
\begin{equation}
\label{vrh}
\int_\Omega v\,dx=b,\qquad
\supp v\subset B(y,Cr),\qquad
\int_\Omega v\cL_0v\,dx\le Cb^2\frac h{r^2}.
\end{equation}
\end{lemma}

\begin{proof}
Fix $y\in\partial\Omega$. After a translation and rotation, we may 
assume $y=0$ and write the boundary near $y$ as
$x_3=\gamma(\xi)$, where $\xi=(x_1,x_2)\in\R^2$,
$\gamma(0)=0$, and $\Omega$ lies above the graph of $\gamma.$
Since $\partial\Omega$ is compact and $C^2$, we may choose
$r_0>0$, independent of $y$, such that
\[
\Phi(\xi,t):=(\xi,\gamma(\xi)+t)\in\Omega
\quad\mbox{if } |\xi|<2r_0 \mbox{ and }  0<t<2r_0,
\]
and $\sup_{|\xi|\le2r_0}|\nabla\gamma(\xi)|\le C$
for a constant $C$ independent of $y$.
Moreover, $\Phi(\xi,0)\in\partial\Omega$,
$\det D\Phi=1$, and $|\Phi(\xi,t)|\le C(|\xi|+t)$.

Choose fixed nonnegative functions
$\eta\in C^\infty_{\rm c}(B_{\R^2}(0,1))$ and
$\zeta\in C^\infty_{\rm c} ((1/4,3/4))$, each with unit integral. Define
\[
v_0(\Phi(\xi,t))
=\frac b{r^2h}\eta(\xi/r)\zeta(t/h),
\]
and extend $v_0$ by zero elsewhere. Then
$v_0\in C^2_{\rm c} (\Omega)$, $v_0\ge0$, and
\[
\int_\Omega v_0\,dx=b,\qquad
\supp v_0\subset B(y,Cr),\qquad
\|v_0\|_{L^\infty(\Omega)}\le\frac{Cb}{r^2h},
\]
where the constant $C>0$ depends only on $\Omega$ and the fixed functions $\eta$ and $\zeta$.

Let $d(x)=\dist(x,\partial\Omega)$. Since
$d(\Phi(\xi,t))\le t$, 
we have $\supp v_0\subset\{d\le3h/4\}$. A sufficiently small
mollification of $v_0$ with a nonnegative kernel of unit integral
gives $v\in C^\infty_{\rm c} (\Omega)$ such that
\[
v\ge0,\qquad
\int_\Omega v\,dx=b,\qquad
\|v\|_{L^\infty(\Omega)}\le\frac{Cb}{r^2h}.
\]
After increasing $C$,
we obtain
$ E:=\supp v\subset\{d<h\}\cap B(y,Cr).$

Decrease $r_0$, if necessary, so that Lemma~3.1 (2) in
\cite{KhanLiFirst} applies with $\ell=h$.
Let $T_E\in H_0^1(\Omega)$ solve
$-\Delta T_E=\1_E$ in $\Omega$.
That lemma gives $0\le T_E\le h^2$. Since
\[
(-\Delta+\sigma)\bigl(\|v\|_{L^\infty(\Omega)}T_E\bigr)
\ge\|v\|_{L^\infty(\Omega)}\1_E\ge v,
\]
the comparison principle yields
\[
0\le\cL_0v
\le\|v\|_{L^\infty(\Omega)}T_E
\le\frac{Cbh}{r^2}.
\]
Consequently,
\[
\int_\Omega v\cL_0v\,dx
\le\frac{Cbh}{r^2}\int_\Omega v\,dx
=Cb^2\frac h{r^2}.
\]
This proves \eqref{vrh}.
\end{proof}

\begin{lemma}
\label{l:mu}
Let $\mu \in \mathcal M(\overline{\Omega})$  be such that $ \mu \ge 0$ on $\overline{\Omega}$, 
$\mbox{supp}\,(\mu) \subseteq \partial\Omega$, and $\| \mu \| = M_{\rm e}.$
\begin{compactenum}
\item[{\rm (1)}]
If $\{u_k\}\subset\cE$ and
$u_k\,dx\stackrel{*}{\rightharpoonup}d\mu$ on
$\overline\Omega$, then
\[
0\le \limsup_{k\to\infty}J[u_k]
\le\int_{\partial\Omega}f(g)\,d\mu =:\ell_{\rm max}.
\]

\item[{\rm (2)}]
Let $\ell$ satisfy $0 \le \ell \le \ell_{\rm max}.$
There exists a sequence
$\{u_k\}\subset\cE\cap C^\infty_{\rm c}(\Omega)$ such that
\begin{equation}
\label{mukmu}
u_k\,dx\stackrel{*}{\rightharpoonup}d\mu
\quad\mbox{on }\overline\Omega \qquad \mbox{and}
\qquad J[u_k]\to \ell.
\end{equation}
If $\ell = \ell_{\rm max},$
then the sequence
can also be chosen so that $ \| u_k \|_{H^{-1}(\Omega)} \to 0.$
\end{compactenum}
\end{lemma}

\begin{proof}
{\rm (1)}
Since $0\le\rho_{u_k}\le\rho_0$,  $f$ is increasing, $f(\rho_0)\in C(\overline\Omega)$, and
$\rho_0=g$ on $\partial\Omega$, we get
\[
0\le J[u_k]=\int_\Omega u_k f(\rho_{u_k}) \,dx\le\int_\Omega u_kf(\rho_0)\,dx
\to\int_{\overline \Omega} f(\rho_0)\, d\mu = \int_{\partial\Omega}f(g)\,d\mu.
\]

{\rm (2)}
We first construct sequences for the two endpoints $\ell=0$ and $\ell=\ell_{\rm max}$, and then use continuity to obtain the
intermediate values.

{\it Step 1.}
We construct a sequence $\{u_k\} \subset \cE\cap C_{\rm c}^\infty(\Omega)$ that satisfies 
\eqref{mukmu} for $\ell = \ell_{\rm max}$ and $\|u_k\|_{H^{-1}(\Omega)}\to0$.  
For each integer $k\ge 1$, cover $\partial\Omega$ by finitely many relatively open sets 
$U_{kj} \subset \partial\Omega$ of diameter less than $1/k$, and then 
set $E_{k1} = U_{k,1}$ and $E_{k j} = U_{k j}\setminus \bigcup_{i < j } U_{k i} $ $(j \ge 2).$
Discard the sets of zero $\mu$-measure and relabel the
remaining sets as $E_{k1},\ldots,E_{kN_k}$.
These are pairwise disjoint Borel sets of positive
$\mu$-measure and diameter less than $1/k$, and they
cover $\partial\Omega$ up to a set of zero $\mu$-measure. Choose $y_{kj}\in E_{kj}$ and place
the mass $\mu(E_{kj})$ at $y_{kj}$. Define

\[
\mu_k=\sum_{j=1}^{N_k}b_{kj}\delta_{y_{kj}},
\qquad
y_{kj}\in E_{kj} \subset \partial\Omega,\qquad b_{kj}= \mu(E_{kj}) >0 \qquad \forall k \ge 1.
\]
We verify for each $k$ that $\sum_{j=1}^{N_k}b_{kj}=M_{\rm e}$, $\mu_k \in \mathcal M(\overline{\Omega}),$
$\mbox{supp}\,(\mu_k) \subseteq \partial \Omega$, and $\| \mu_k \| = M_{\rm e}$. 
Moreover, for every $\varphi\in C(\overline\Omega)$,
\begin{align*}
\left|
\int_{\overline{\Omega}}\varphi\,d\mu_k
-\int_{\overline{\Omega}}\varphi\,d\mu
\right|& = \left|
\int_{\partial\Omega}\varphi\,d\mu_k
-\int_{\partial\Omega}\varphi\,d\mu
\right|
\\
&=
\left|
\sum_{j=1}^{N_k}
\int_{E_{kj}}
\bigl[\varphi(y_{kj})-\varphi(y)\bigr]\,d\mu(y)
\right|\\
&\le M_{\rm e}
\sup_{{x,y\in\partial\Omega, |x-y|\le1/k}}
|\varphi(x)-\varphi(y)|
\to 0,
\end{align*}
by the uniform continuity of $\varphi$.
Thus $\mu_k\stackrel{*}{\rightharpoonup}\mu$ on
$\overline\Omega$.

Let $r_0,C$ be as in Lemma~\ref{l:construction}.
Choose $\{r_k\}\subset(0,\min\{r_0,1\})$ with
$r_k\to0$, and set $h_k=r_k^3$.
The lemma gives nonnegative
$v_{kj}\in C^\infty_c(\Omega)$ satisfying
\begin{equation}
\label{vkjvkjvkj}
\int_\Omega v_{kj}\,dx=b_{kj},\qquad
\supp v_{kj}\subset B(y_{kj},Cr_k),\qquad
\int_\Omega v_{kj}\cL_0v_{kj}\,dx\le Cb_{kj}^2r_k.
\end{equation}
Set $u_k=\sum_{j=1}^{N_k}v_{kj}$. Then
$u_k\in\cE\cap C^\infty_c(\Omega)$.
For $\varphi\in C(\overline\Omega)$, we have
\[
\begin{aligned}
\left|\int_{\overline{\Omega}} u_k\varphi\,dx
-\int_{ \overline{\Omega}}\varphi\,d\mu_k\right|
&\le\sum_{j=1}^{N_k}
\int_\Omega v_{kj}|\varphi(x)-\varphi(y_{kj})|\,dx\\
&\le M_{\rm e}
\sup_{{x,y\in\overline\Omega, |x-y|\le Cr_k}}
|\varphi(x)-\varphi(y)|
\\
& \to 0 \qquad\mbox{as } k \to \infty.
\end{aligned}
\]
Since $\mu_k\stackrel{*}{\rightharpoonup}\mu$, the preceding
estimate implies
$u_k\,dx\stackrel{*}{\rightharpoonup}d\mu$ on $\overline\Omega$.

With the norm $ \| \cdot \|_\sigma $ defined in Lemma~\ref{l:L0}, we have
$\|u_k \|_\sigma = \| \sum_{j=1}^{N_k} v_{kj} \|_\sigma \le \sum_{j=1}^{N_k} \| v_{kj} \|_{\sigma}$. Hence, 
by Lemma~\ref{l:L0} and \eqref{vkjvkjvkj}, 
\begin{align*}
    \int_\Omega u_k\cL_0u_k\,dx = \| u_k \|^2_\sigma 
\le  \left[ \sum_{j=1}^{N_k} \left( \int_{\Omega } v_{kj} \cL_0 v_{kj} \, dx \right)^{1/2} \right]^2
\le Cr_k \left( \sum_{j=1}^{N_k} b_{kj} \right)^2 
=CM_{\rm e}^2r_k
\to 0.
\end{align*}
Since by Lemma~\ref{l:L0}, the norms $\| \cdot \|_\sigma $ and $\| \cdot \|_{H^{-1}(\Omega)}$
are equivalent on $L^2(\Omega)$, we have $ \|u_k \|_{H^{-1}(\Omega)} \to 0.$
Using \eqref{Jvrho0} with $u_k$ replacing $v$, we obtain
\[
0\le\int_\Omega u_kf(\rho_0)\,dx-J[u_k]
\le f(g_{\rm max}) \| u_k \|_\sigma^2
\to 0.
\]
Note that
\[
\int_\Omega u_kf(\rho_0)\,dx =
\int_{\overline{\Omega}} u_k f(\rho_0)\, dx
\to\int_{\overline{\Omega}} f(\rho_0) \, d\mu = \int_{\partial\Omega}f(g)\,d\mu = \ell_{\rm max}.
\]
The limit $J[u_k]\to \ell_{\rm max}$ follows from the above two limits.

{\it Step 2.}
We construct a sequence
$\{w_k\}\subset\cE\cap C^\infty_{\rm c}(\Omega)$ with the same limiting measure $\mu$
such that  $J[w_k]\to0$.
From Step~1, choose distinct
$x_{kj}\in\Omega$ with $|x_{kj}-y_{kj}|<1/(2k)$.
Choose $r_{kj}>0$ such that
\begin{align*}
&B(x_{kj},2r_{kj})\subset\subset\Omega\cap B(y_{kj},1/k),\\
&
\sum_{j=1}^{N_k}r_{kj}<\frac1k,\\
&B(x_{ki},2r_{ki})\cap B(x_{kj},2r_{kj})=\varnothing
\quad(i\ne j).
\end{align*}
Fix a nonnegative $\eta\in C^\infty_{\rm c}(B(0,1))$ of integral
one, and define
\[
w_{kj}(x)=\frac{b_{kj}}{r_{kj}^3}
\eta\left(\frac{x-x_{kj}}{r_{kj}}\right)
\qquad \mbox{and} \qquad 
w_k=\sum_{j=1}^{N_k}w_{kj}.
\]
Then $w_k\in\cE\cap C^\infty_{\rm c} (\Omega)$ and, as in Step~1,
$w_k\,dx\stackrel{*}{\rightharpoonup}d \mu$.

Choose $\chi\in C^\infty_{\rm c} (B(0,2))$ with
$0\le\chi\le1$ and $\chi=1$ on $B(0,1)$, and set
\[
\chi_{kj}(x)=\chi\left(\frac{x-x_{kj}}{r_{kj}}\right)
\qquad \mbox{and} \qquad 
\chi_k=\sum_{j=1}^{N_k}\chi_{kj}.
\]
Then $0\le\chi_k\le1$, $\chi_k=1$ on $\supp w_k$, and
$
\|\Delta\chi_{kj}\|_{L^1(\Omega)}
=r_{kj}\|\Delta\chi\|_{L^1(\R^3)}.
$
Multiplying the equation for $\rho_{w_k}$ by $\chi_k$
and integrating by parts, we obtain
\[
\begin{aligned}
0& \le  J[w_k]
=\int_\Omega w_k f(\rho_{w_k})\, dx 
=\int_\Omega w_kf(\rho_{w_k})\chi_k\,dx
=\int_\Omega\rho_{w_k}\Delta\chi_k\,dx
-\sigma\int_\Omega\rho_{w_k}\chi_k\,dx\\
&
\le g_{\rm max}
\sum_{j=1}^{N_k}\|\Delta\chi_{kj}\|_{L^1(\Omega)}
\le Cg_{\rm max}\sum_{j=1}^{N_k}r_{kj}
<\frac{Cg_{\rm max}}k \to 0 \qquad \mbox{as } k\to \infty.
\end{aligned}
\]

{\it Step 3.}
Let $0 < \ell < \ell_{\rm max}.$
By Steps~1 and 2, after discarding finitely many terms
and relabeling, we have
$J[w_k]<\ell<J[u_k]$ for all $k \ge 1.$  
For each $k\ge 1$, 
Proposition~\ref{p:Juv} and the Intermediate Value Theorem
applied to the function $t \mapsto J[(1-t) w_k + t u_k]$
$(0 \le t \le 1)$ 
give $t_k\in(0,1)$ such that
$
J[(1-t_k)w_k+t_ku_k]=\ell.
$
Since both $\{u_k\,dx\}$ and $\{w_k\,dx\}$ converge weakly-$\ast$
to $d \mu$ and $0<t_k<1$, we also have
\[
\bigl((1-t_k)w_k+t_ku_k\bigr)\,dx
\stackrel{*}{\rightharpoonup}d\mu.
\]
Together with Steps~1 and 2, this proves {\rm (2)}.
\end{proof}

\bibliography{Reaction_Flux_Landscape_References}
\bibliographystyle{plain}

\end{document}